\documentclass[12pt]{article}

\usepackage{amssymb,amsmath}
\usepackage{amsfonts}
\usepackage[T1]{fontenc}
\usepackage{xcolor}
\usepackage{ulem}

\newtheorem{Theorem}{Theorem}[section]
\newtheorem{Lemma}[Theorem]{Lemma}
\newtheorem{Corollary}[Theorem]{Corollary}
\newtheorem{Proposition}[Theorem]{Proposition}
\newtheorem{Note}[Theorem]{Remark}
\newtheorem{Example}[Theorem]{Example}
\newtheorem{Conjecture}[Theorem]{Conjecture}
\def\ep{\hfill{\vbox to 7pt{\hbox to 7pt{\vrule height 7pt width 7pt}}}}

\def\A{(A,\omega)}
\def\train{x^n+\gamma_1\omega(x)x^{n-1}+\ldots+\gamma_{n-1}\omega(x)^{n-1}x}
\renewcommand*\descriptionlabel[1]{\hspace\labelsep
	\normalfont\rmfamily #1}

\numberwithin{equation}{section}

\newcommand{\Q}{{\mathbb Q}}
\newcommand{\N}{{\mathbb N}}
\newcommand{\Z}{{\mathbb Z}}
\newcommand{\R}{{\mathbb R}}
\newcommand{\Complex}{{\mathbb C}}

\def\l{\lambda}
\def\esp{\hspace*{0,5cm}}
\def\ep{\hfill{\vbox to 7pt{\hbox to 7pt{\vrule height 7pt width 7pt}}}}
\def\pr{{\it Proof. }\ }
\addtocounter{section}{-1}

\newcommand{\AdOper}{Ad}
\newcommand{\al}{\alpha}
\newcommand{\alg}{\mathop{\mathrm{alg}}\nolimits}
\def\bc{\begin{center}}
\newcommand{\be}{\beta}
\newcommand{\charact}{\mathop{\mathrm{char}}\nolimits}
\newcommand{\de}{\delta}
	\def\ec{\end{center}}
\newcommand{\goth}{\mathfrak}
\newcommand{\Identity}{Id}
\newcommand{\iso}{\,\tilde \rightarrow\,}
\newcommand{\g}{\gamma}
\newcommand{\la}{\lambda}
\newcommand{\La}{\Lambda}
\newcommand{\om}{\omega}
\newcommand{\va}{\varphi}
\DeclareMathOperator*{\sd}{\times}
\newcommand{\sig}{\sigma}

\def\lan{\langle }
\def\ran{\rangle }

\begin{document}
	\centerline{\Large\bf  Bernstein algebras that are algebraic}
	\vspace*{0,4cm}
	\centerline{\Large\bf  and the Kurosh problem}
	
	\vspace*{0,2cm}
	
	\vspace*{1cm} \centerline{\bf   Dmitri PIONTKOVSKI\footnote{Corresponding author}}
	\vspace*{0,2cm}
	\centerline{\small
		HSE University}
	\centerline{\small 20 Myasnitskaya Str.
		Moscow 101000 Russia}
	\centerline{\small  E-mail:  dpiontkovski@hse.ru}
	\vspace*{0,7cm}
	\centerline{\bf  Fouad ZITAN}
	\vspace*{0,2cm}
	\centerline{\small
		Laboratory of Algebra and its Applications,
		Faculty of
		 Sciences,}
	\centerline{\small
		Abdelmalek Essaadi University, BP $2121$, T\'etouan, Morocco}
	\centerline{\small  E-mail: fzitan@uae.ac.ma}
	
	\vspace*{1cm}
	\centerline{\it To the memory of Serge Bernstein (1880--1968)}
	\centerline{\it whose works were the source of Bernstein algebras}
	\vspace*{0,5cm}
	\centerline{\it In honor of Professor Yuri I. Lyubich}
	\centerline{\it on the occasion of his 90th birthday}
	
	\vspace*{0,5cm}
	\baselineskip0.5cm
	\centerline {{{\bf Abstract}}} \vskip 5mm
	{\small \esp 
	We study the class of 
	Bernstein algebras that are algebraic, in the sense that each element generates a finite-dimensional subalgebra. Every Bernstein algebra has a maximal algebraic ideal, and the quotient algebra is a zero-multiplication algebra. Several equivalent conditions for a Bernstein algebra to be algebraic are given. In particular, known characterizations of train Bernstein algebras in terms of nilpotency are generalized to the case of locally train algebras. Along the way, we show that  if a Banach Bernstein algebra is algebraic (respectively, locally train), then it is of bounded degree (respectively, train).
	
Then we investigate the Kurosh problem for Bernstein algebras: 
		whether a finitely generated Bernstein algebra which is algebraic of bounded degree is finite-dimensional. This problem turns out to have a closed link with a question about associative algebras. In  particular, when the barideal  is nil, the Kurosh problem
asks whether a finitely generated Bernstein-train algebra is finite-dimensional.  We prove that the answer is positive for some specific cases and for low degrees, and construct counter-examples in the general case. 
		
      By results of Yagzhev, the Jacobian conjecture is equivalent to a certain statement about  Engel and nilpotence identities of multioperator algebras. We show that 
    the generalized Jacobian conjecture for quadratic mappings
    holds for Bernstein algebras. 
		}\\
	
		

	\baselineskip0.6cm
	
	\bigskip
	
	\noindent {\cal \emph{AMS}} \emph{Classification}: 17D92
	
	\noindent {\cal \emph{Keywords}}:  Algebraic element; Banach algebra; Bernstein algebra; Finitely generated algebra; Jordan algebra; Kurosh problem; Locally nilpotent
	algebra; Train algebra.

	
	
	\section{\bf Introduction}
	
	\baselineskip0.6cm
		\esp Let $Q$ be a quadratic operator on a vector space $A$ (over an infinite field $K$ with $\charact K \ne 2$), so that the polarization
	of $Q$ defines a commutative (non-associative) multiplication on $A$ such that $x^2 = Q(x)$. Suppose that $Q$ maps some hyperplane $H$ to itself and satisfies a stationary condition
	$Q^2(s) = Q(s) $ for all $s\in H$. Then the algebra $A$ is called a Bernstein algebra.

	Bernstein algebras 
	give 
	an
	algebraic formulation of the population genetics problem of classifying the stationary
	evolution operators  $Q$ posed by Bernstein (see \cite{Berns1, Berns2}).
	Let $p = (p_1, \dots, p_n) \in \R^n$ be the distribution of frequencies of $n$ genes in a population of some species, 
	and let $Q(p)$ be the distribution of the same genes in the next generation after linkage. In 1908, Hardy 
	and Weinberg~\cite{Hardy, Weinberg} independently deduced from the Mendelian laws that under some natural assumptions, the distribution of genes becomes stable after the second (but not the first) 	linkage. Nearly a century ago, Bernstein reformulated the Hardy--Weinberg law 
	as the elegant stationary equation
	$Q^2(s) = Q(s) $ for the quadratic operator $Q$ and all elements $s$ of the 
	standard simplex $S=\{ (x_1, \dots ,x_n )| x_1 +\dots +x_n=1, x_i \ge 0 \}$.
	He also posted a challenging problem to classify such genetic systems; Bernstein himself had solved it in low dimensions. 
	 About half a century later, Lyubich (see~\cite{Ly} and references therein; see also~\cite{Ho}) introduced Bernstein algebras and gave a classification for all dimensions under mild additional conditions. 
	In the Lyubich approach, the problem of the description of such an operator $Q$ in $\R^n$ corresponds to the one of  the description of $n$-dimensional Bernstein algebras. 
	Since the number $n$ of genes may be large enough, the limit case of infinite-dimensional Bernstein algebras is also of interest.

	Bernstein
	algebras have been the subject of several researches (see, for
	instance, \cite{Gonz}, \cite{Ly}, \cite{Mi}), and there is at
	present a substantial bibliography on the subject.
	Our knowledge of infinite-dimensional Bernstein
	algebras is recent and began with finitely generated Bernstein
	algebras. Peresi \cite{Pe} and Krapivin \cite{Kra} have shown
	independently that the barideal of a finitely generated nuclear algebra is
	nilpotent and so finite-dimensional  \cite{Su}. Later, Boudi and Zitan \cite{B, Zitan} 
	undertook a systematic study of Bernstein
	algebras satisfying chain conditions on ideals.
	
	 Another kind of non-associative algebras arising in  the symbolism of genetics are the so-called train algebras, which
	were introduced by Etherington \cite{Eth1, Eth2}. The class of Bernstein algebras which are train algebras were studied by several authors (see, for instance, \cite{GarciaGonzalez, Ouatt, Wa1}), but mostly restricted to the finite-dimensional context.\\
	
	In this paper, we study  algebraic and nilpotent elements in Bernstein algebras. 
	We begin with a discussion 
	on 
	various versions of algebraicity and nilpotency properties in
	Bernstein algebras, and study in details the case of a low degree of algebraicity.   
	We say that a Bernstein algebra is algebraic if each  element $a$ is algebraic, that is,
	the subalgebra $\alg(a)$ generated by it is finite-dimensional. In a Bernstein algebra, every algebraic element satisfies an equation over the ground field 
	\begin{equation}
	\label{eq:intro:algebraic_element}
	a^n - \gamma_{n-1} a^{n-1}-\dots - \gamma_{1} a = 0,
	\end{equation}
	where $a,a^2, \dots, a^{n} = a a^{n-1}$ are the principal powers of $a$.
    In a sense, each Bernstein algebra is closed to being algebraic, as the next theorem shows. 
		\begin{Theorem}[Theorem~\ref{th:alg_ideal}]
	\label{th:intro:alg_ideal}
	Any Bernstein algebra $A$ contains a maximal algebraic ideal ${\cal I}(A)$. Furthermore,  
	${\cal I}(A)$ contains $A^2$, so that $A/{\cal I}(A)$ is a zero-multiplication algebra. 
	\end{Theorem}	
	
	We show that a 	Bernstein algebra is algebraic (of bounded degree) if and only if some corresponding
	multiplication operators are locally algebraic (algebraic of bounded degrees collectively). Many results in this direction are given, involving linear operators.
		Particularly, it is proved that every algebraic Bernstein algebra of degree $\leq 2$ is a Jordan algebra, and counter-examples for degree $\geq3$ are provided. 
	
	 In general, a Bernstein algebra cannot be nil since the elements of $H$ are not nilpotent. Still, there is an ideal in such an algebra which is nil in the important cases.
	Any Bernstein algebra $A$ admits a homomorphism $\omega$ onto $K$ (that is, $A$ is {\em baric}) such that $\omega(s)=1$ for all $s\in H$. The kernel $N = \ker \omega$ is called the {\em barideal} of $A$ and plays an important role in the theory. For each $s\in H$, we have $N=H-s$.

		Recall that a baric algebra is train if all  elements of $H$ satisfy the same polynomial equation
		of the form~(\ref{eq:intro:algebraic_element}). 
		In the relaxed version, a baric  algebra is called locally train~\cite{Ouatt1} if
	each singly-generated subalgebra is train. 
		Algebraic Bernstein algebras (of bounded or unbounded degrees, respectively) are natural
	generalizations of Bernstein train  and locally train algebras. 
	For finite-dimensional Bernstein algebras,
	both train and locally train conditions  are equivalent to the fact that $N$ is a nil-subalgebra.
	Moreover, in this case it is equivalent to say that the  Engel identity
	$$ 
	x(x(\dots (xy)\dots )) = 0
	$$
	of some degree holds in $N$, that is, $N$ is an Engel subalgebra. We 
	establish 
	\begin{Theorem}
	\label{th:intro:train}
	 With the notation above, the following conditions are equivalent:
	
	(i) the Bernstein algebra $A$ is (locally) train;
	
	(ii) for each $x\in N$, the Engel identity of some degree holds for all $y\in N$ 
	(respectively, $N$ is Engel);
	
	(iii) $N$ is a nil-algebra (of bounded nil-index).
	\end{Theorem}
	
	Furthermore, we use a connection with the Jacobian conjecture in Theorem~\ref{th:Bernstein_Engel_Yagzhev} to establish that the equivalence $(ii)\Longleftrightarrow (iii)$ above holds for arbitrary subalgebras of a Bernstein algebra in place of $N$. 
	
	In particular, all the above conditions are equivalent in the case of Banach Bernstein algebras: this follows from  
\begin{Theorem}[Theorems~\ref{th:banach1} and~\ref{th:banach2}]
    \label{th:intro:banach}
    Let $A$ be a Banach Bernstein algebra.
    
    1. $A$ is algebraic if and only if it is algebraic of bounded degree.
    
    2. $A$ is train if and only if it is locally train.
\end{Theorem}
	
	Next,
	we examine the following three interesting problems in the class of Bernstein algebras.
	
	{\em The Kurosh problem.} Whether a finitely generated Bernstein algebra is necessarily finite-dimensional provided that each element is algebraic? 
	The Kurosh problem has been intensely studied in various classes of algebras, as associative, Jordan, alternative and Lie algebras. Thus, it is quite legitimate to investigate
	the counterpart of Bernstein algebras with regard to this problem. When the barideal $N$ is nil, the Kurosh problem asks whether a finitely generated Bernstein-train algebra is finite-dimensional. We will see that the Kurosh problem for Bernstein
	algebras is equivalent to a weak form of the Kurosh problem for associative algebras. We solve positively
	some special cases, and construct a family of counter-examples in the general case.
	
	{\em The Burnside-type problem}. In the class of associative algebras, the standard 
	Burnside-type problem is whether each finitely-generated nil-algebra is nilpotent. 
	For finitely generated Bernstein algebras, one can naturally ask: 
	If each element of the algebra generates a train subalgebra, is the algebra itself  necessarily train? 
	
	We solve both the Kurosh and the Burnside problems for Bernstein algebras. The  answers
	are positive for algebras of  low degrees of algebraicity, but negative in general.
	
    	\begin{Theorem}[Theorem~\ref{th:Kurosh}]
    	\label{th:intro-Kurosh}
		All finitely generated Bernstein algebras which are
		algebraic of bounded degree $\leq 2$ (in particular,   Bernstein-train algebras of rank at most $3$) are  finite-dimensional.
		
		In contrast, there exists  an infinite-dimensional finitely generated Bernstein train algebra of rank $4$. It follows that for all $d\ge 3$ and $n\geq 2$, there exists an  infinite-dimensional  finitely generated Bernstein algebra with $n$ generators, which is
		algebraic of bounded degree $d$. 
	\end{Theorem}
	
	In fact, the algebra in the second claim of the theorem can be chosen to be a train algebra of rank $d+1$, Cf. Remark~\ref{rem:Kurosh_example_n_k}.
	
	{\em The Jacobian problem}. The famous Jacobian problem asks whether each polynomial 
	endomorphism of $\Complex^n$ with constant Jacobian determinant has a polynomial inverse.
	Due to results of Yagzhev (\cite{Yagzhev}; see also a survey in~\cite{Belov_etc}), the Jacobian problem 
	is equivalent to an implication of the form: if a (multi-operator) algebra is Engel and satisfies a system of Capelli identities of degree $n+1$, then it is weakly nilpotent. 
	The {\em generalized Jacobian conjecture for quadratic mappings}~\cite{Belov_etc} claims that binary Engel algebras are 
	{\em weakly}
	nilpotent: this is a form of the Jacobian conjecture for quadratic polynomials and infinite numbers of variables. Whereas the Jacobian conjecture holds for quadratic polynomials~\cite{Wang}, its infinite-dimensional version is still open. We prove  it under the additional assumption that the algebra is Bernstein. 
	\begin{Theorem}[Corollary~\ref{cor:jacobian-bernstein}]
\label{th:intro-jacobian-bernstein}
	Suppose that $A$ is a subalgebra of a Bernstein algebra. Then the  generalized Jacobian conjecture for quadratic mappings holds for $A$.  
\end{Theorem}

The same theorem is proved for commutative algebras satisfying the identity $(x^2)^2=0$, see Theorem~\ref{th:Bernstein_Engel_Yagzhev}.\\

	The paper is organized as follows. In Section~\ref{sec:prelim} we recall the basic theory of Bernstein algebras. 
	We use the following definition which is equivalent to the one discussed above.
	A baric  commutative  algebra $A$
	(over an infinite field $K$ of characteristic different from 2
	and 3) is called Bernstein if 
	the elements of the subset $H = \{ x \;|\; \omega(x) =1\}$ satisfy the equation 
	$Q(Q(x)) = Q(x)$, where $Q(x) = x^2$. Equivalently, such an algebra satisfies the identity $(x^2)^2=\omega(x)^2 x^2$. Then we recall the definitions of the Peirce decomposition
	and the Lyubich ideal in Bernstein algebras as well as some descriptions of special types of such algebras.
	
	
	In Section~\ref{sec:alg_Bernstein_algebras}, we discuss algebraic Bernstein algebras. 
	Among others, we show that a Bernstein algebra $A$ is algebraic (of bounded degree) if and only if either the elements of $H$ or $N$ are algebraic (respectively, of bounded degree). We also prove  Theorem~\ref{th:intro:alg_ideal}. 

	In Section~\ref{sec:singly_generated}, we give a description of a Bernstein algebra generated by a single element. In particular, we show in Remark~\ref{rem:single-generated-free} that there exists a free single-generated Bernstein algebra (since Bernstein algebras do not form a variety, there is no general concept of free  Bernstein algebras). 
	We  utilize these results in Section~\ref{sec:locally_train}, where we examine train and locally train Bernstein algebras.  
	In particular, we prove here Theorem~\ref{th:intro:train}.
		Other equivalent conditions for train and locally train algebras are also given.  
	
	Section~\ref{sec:banach} deals with Banach Bernstein algebras, which were introduced by Lyubich~\cite{Ly2}. We apply the results of
	Cabrera and Rodr\'iguez~\cite{Rod} to establish Theorem~\ref{th:intro:banach}.
 Section~\ref{sec:alg_low_degree} is devoted to Bernstein algebras which are algebraic of low degrees. In particular, Corollary~\ref{cor:former4.4} shows that Jordan Bernstein algebras are exactly those algebraic algebras of degree at most 3.
In Sections~\ref{sec:Kurosh} and~\ref{sec:Jacobian}, we discuss the Kurosh and the Jacobian problems.  In particular, we prove here Theorems~\ref{th:intro-Kurosh}
and~\ref{th:intro-jacobian-bernstein}.\\
Various examples are presented throughout the article to serve as motivation and illustration for our results.


	
	\section{Preliminaries}
	
	\label{sec:prelim}
	
	\esp In order to keep the paper reasonably self-contained, in this section we
	summarize the basic notions that will be used in
	this work. Other notions and properties will be given where
	appropriate.\\  	
	Let $K$ be an infinite field of characteristic different from 2
	or 3,  and let $A$ be an algebra over $K$, not necessarily
	associative or finite-dimensional. If $A$ has a non-zero algebra 
	homomorphism $\omega: A \rightarrow K$, then the ordered pair $(A, \omega)$
	is called {\it a baric algebra} and $\omega$ is its {\it weight
		function}. For each $e\in A$ with $\omega(e)\neq 0$, we have
	$A=Ke\oplus N$, where $N=\ker(\omega)$ is an ideal of $A$, called
	the {\it barideal} of $A$. A baric ideal of A is an ideal $I$ of A with $I\subset N$. In this case, the quotient algebra $A/I$
	is a baric algebra with weight function $\overline{\omega}$ defined by $\overline{\omega}(x+I)=\omega(x)$.
	
	{\it A Bernstein algebra} is a commutative baric algebra $(A,\omega)$
	satisfying the identity $(x^2)^2=\omega(x)^2 x^2$. In the notation of the Introduction, 
	$H = \{x\in A \;|\; \omega (x)=1\}$
	and $Q(x) = x^2$. If $x\in H$, 
	then $e=x^2$ is a nontrivial idempotent of $A$ that gives rise to the
	Peirce decomposition $A=Ke\oplus U\oplus V$, where $N=\ker(\omega)=U\oplus V$,
	and
	\begin{equation}
		U=\{u\in A\; |\; 2eu=u\},\; \; V=\{v\in A\; |\; ev=0\}.
	\end{equation}
	We now collect some well-known results about Bernstein algebras (see,
	for instance, \cite{Gonz}, \cite{Ly}, \cite{Mi}, \cite{Wo}). For
	every idempotent $e$ of the Bernstein algebra $\A$, the Peirce
	components multiply according to
	\begin{equation}
	\label{eq:Bernstein_UV}
		U^2\subseteq V, \; UV\subseteq U, \; V^2\subseteq U, \; UV^2=0.
	\end{equation}
	We also have the identities for $u\in U,\, v\in V$ and $x, y\in
	N$:
	\begin{equation}
	\label{eq:Bernstein_identities=0}
	u^3=u(uv)=u^2(uv)=(uv)^2=u^2v^2=0.
	\end{equation}
	
	In addition, we have \cite{Mi}:
	\begin{equation}U^iV^j=0 \mbox{ for all } i\geq 1 \mbox{ and } j\geq 2\end{equation}
	
	The set of idempotents of $A$ is given by $I(A)=\{e+u+u^2 \;
	| \; u\in U\}$.\\ If $e'=e+u_0+u_0^2$ is another idempotent, the
	attached Peirce components are
	\begin{equation} U'=\{u+2u_0u \;|\; u\in U\}, \;\;
		V'=\{v-2(u_0+u_0^2)v \; | \; v\in V\}.\end{equation}
	 When $A$ is finite-dimensional,
the pair of integers $(1 + r, d)$, where $r = \dim U$ and $d = \dim V$, is an invariant called the {\it type} of $A$.
	A Bernstein algebra $A=  Ke\oplus U\oplus  V$ is not unital unless in the trivial case $\dim A= 1$, and
	cannot be associative except when $U=0$. However, Bernstein algebras may be {\it power-associative},
	that is, if each element generates an associative subalgebra.\\
	
	Recall that a commutative algebra $A$ is a {\it Jordan algebra} if the identity $x(x^2y)= x^2(xy)$ holds in $A$. Bernstein Jordan algebras play a fundamental role in the theory of Bernstein algebras. It is well known
	that the following four conditions are equivalent for a Bernstein algebra $A= Ke\oplus U\oplus V$ (see, for instance, \cite{Wa}, \cite[Corollary 3.5.17]{Ly}):\\
	(a) $A$ is a Jordan algebra.\\
	(b) $A$ is power-associative.\\
	(c) $x^3=\omega(x)x^2$ for all $x\in  A$.\\
	(d) $V^2=0$ and $(uv)v= 0$ for all $u\in  U$ and $v\in  V$.\\
	
	Therefore, the elements of the barideal $N= \ker(\omega)=U\oplus V$ in a Bernstein Jordan algebra $(A, \omega)$ satisfy
	$x^3=0$ and so the Jacobi identity
	$(xy)z + (yz)x + (zx)y = 0$.\\
	
		The subspace $L(A) = ann_U(U) =\{u\in U \; |\; uU=0\}$ is an ideal of $A$,
	called {\it the Lyubich ideal} (or {\it the Jordan ideal}). It is independent of the
	chosen idempotent $e$, and satisfies: \begin{equation}L(A) (U\oplus U^2)=0, \;
		V^2\subseteq L(A) ,\;  v(vu)\in L(A)  \mbox{ for all } u\in U \mbox{ and } v\in V \end{equation} (see \cite{Gonz} and \cite[Theorem 3.4.19]{Ly}). This
	ideal $L(A) $ plays a decisive role in the connection between Bernstein
	and Jordan algebras, since the factor algebra $A/L(A) $ is a
	Bernstein Jordan algebra (see, for instance, \cite{Gonz} and
	\cite[Corollaire 3.3]{Mi}). 
	
	Remember that if $I$ is an ideal of a Bernstein
	algebra $A$, then by \cite[Proposition 2.1]{Koul}, either $I=(I\cap
	U)\oplus (I\cap V)$ or $I=Ke\oplus U\oplus V'$, where $V'$ is a
	subspace of $A$ with $U^2\subseteq V'\subseteq V$. 
	
	A Bernstein algebra $A$ is called {\it
		nuclear} if $A^2=A$, or equivalently, $U^2=V$. In this case, we
	have $N=\ker(\omega)=U\oplus  U^2$ and $L(A) N=0$. Moreover,
	$A$ is the only ideal which is not
	baric.\\ A Bernstein algebra $A=Ke\oplus U\oplus V$ is said to be
	{\it exceptional} (or {\it exclusive}) if it satisfies the (invariant) condition
	$U^2=0$, or equivalently, $(xy)^2=\omega(xy)xy$ for all $x,y\in A$. It is important to realize that any Bernstein algebra
	$A=Ke\oplus U\oplus V$ contains a nuclear Bernstein subalgebra
	$A^2=Ke\oplus U\oplus U^2$ and an exceptional Bernstein subalgebra
	$B=Ke\oplus L(A) \oplus
	V$, which are the largest ones.\\
	
	In another context, a commutative baric algebra $(A,\omega)$ is called a {\it train algebra of rank $r$}
	if there exist scalars $\gamma_1,\dots,\gamma_{r-1}$ in $K$ such that
	\begin{equation}
	\label{train}
	x^r + \gamma_1\omega(x) x^{r-1} +
		\cdots+\gamma_{r-1}\omega(x)^{r- 1}x=0,\end{equation} for all $x\in
	A$, where $r\geq 2$ is the smallest integer for which such an
	equation holds, and $x^1=x, \dots, x^{k+1}=x^kx$ are the {\it
		principal powers} of $x$. Equation \eqref{train} is called {\it
		the train equation} of $A$, where we have necessarily
	$1+\gamma_1+\cdots+\gamma_{r-1}=0$. Then  $\ker(\omega)$ is the set of nilpotent
	elements. Consider the ordinary polynomial $P(X)= X^r +
	\g_1X^{r-1}+\cdots + \g_{r-1}X\in K[X]$, called {\it the train polynomial}
	of $A$. In a suitable extension of $K$, $P(X)$ splits into linear
	factors $P(X)=X(X-1)(X-\la_1)\cdots(X-\la_{r-2})$, where $\la_0=1,
	\la_1, \dots, \la_{r-2}$ are called  {\it the principal train roots}
	of $A$. In an abuse of notation as in \cite{Wo}, we express
	\eqref{train} in the form
	$x(x-\omega(x))(x-\lambda_1\omega(x))\cdots(x-\lambda_{r-2}\,\omega(x))=0$, which really
	means $\left(L_x-\omega(x)\,\texttt{id}_A\right)\left(L_x-\lambda_1
	\omega(x)\,\texttt{id}_A\right)\cdots 	\left(L_x-\lambda_{r-2}\,\omega(x)\,\texttt{id}_A\right)x=0\,$ for all $x\in
	A$, where $L_x$ indicates the multiplication by $x$ and
	$\texttt{id}_A$ stands for the identity mapping.\\
	When $r=3$,  the train equation
	is $x^3-(1+\g)\omega(x)x^2+\g\omega(x)^2x=0$, with $\g\in K$.
	Thus, Bernstein Jordan algebras are special
	instances of train algebras of rank 3 with $\gamma=0$.\\
	
	A train algebra $(A, \omega)$ of rank 2 satisfies the identity $x^2=\omega(x)x$, or equivalently
	$xy=\frac 12(\omega(x)y+\omega(y)x)\; $ for all $x,y\in A$. Then $A$ is a Bernstein algebra, called  {\it  elementary} (or, {\it  unit}). It is also called a {\it gametic algebra of simple Mendelian inheritance}. In this case, $A$ is a Jordan algebra \cite[page 138]{Eth2}; and more precisely $A$ is a special Jordan algebra \cite[page 133]{Shafer}.\\

	Further information about the algebraic properties
	of Bernstein and train algebras, as well as their possible genetic interpretation can be found in \cite{Ly,Reed, Wo}.\\
	
	\esp We now let $A$ be an arbitrary algebra over $K$. We say that $A$ is {\it nilpotent} ({\it right nilpotent}) if the
	descending chain of ideals (right ideals) defined
	recursively by $A^1= A$ and $A^n=\sum\limits_{i+j=n}
	A^i A^j$ ($A^{\lan 1 \ran}=A$ and $A^{\lan i \ran}=A^{\lan i-1 \ran}A$) ends up in zero.
	Clearly, if $A$ is
	nilpotent, then $A$ is right
	nilpotent. Conversely, if $A$ is commutative, then $A^{2^n}\subseteq A^{\lan n \ran}$ by \cite[Proposition 1]{russe}. Therefore,
	if $A$ is a right nilpotent commutative algebra, it is nilpotent too.
	We define {\it the plenary powers} of $A$  by
	$A^{(1)}=A^2$ and $A^{(n)}=(A^{(n-1)})^2$. The algebra $A$ is
	said to be {\it solvable} when $A^{(n)}=0$ for some $n$, and the smallest such $n$ is 
	{\it the index of solvability} of $A$.
	
	For any subset $X\subset A$, we write $\langle X\rangle $ for the subspace of
	$A$ spanned by $X$, and $\text{alg}(X)$ for the subalgebra of $A$ generated by $X$. The algebra $A$ is {\it finitely
		generated} if it is generated as an algebra by a finite subset of $A$.\\
	An algebra  $A$ is {\it locally nilpotent} if
	every finitely generated subalgebra of $A$  is nilpotent.\\
	
	 Returning to Bernstein algebras, recall that in a Bernstein algebra, the principal powers $N^{\langle i \rangle}$
	are ideals \cite[page 113]{Ly}. Moreover, the barideal $N$ satisfies the equation $(x^2)^2=0$, but is not in
	general nilpotent. However, $N$ is always solvable, since $N^{(3)}= 0$ \cite[Theorem 2.11]{Bernad} (see
	also \cite{Jac}).\\
	
	The following  fact, which is easily checked, will be used 
	in this text. It furnishes  sufficient (but not necessary) conditions for an arbitrary commutative algebra to be an exceptional Bernstein algebra. See the proof in \cite[Theorem 3.4.15 (2)]{Ly}.
	\begin{Lemma} 
		\label{lemma:Lemma1.1}
		Let $A$ be an arbitrary commutative algebra containing a nonzero idempotent $e$ and two subspaces $U$ and $V$ satisfying $$A=Ke\oplus U\oplus V, \;\;eu=\frac 12u \;\,\forall u\in U, \;\;ev=0\;\,\forall v\in V,\;\; UV\subset U,\;\; V^2\subset U,\;\; U^2=0.$$ Then $A$ is an exceptional  Bernstein algebra with Peirce components $U$ and $V$ with respect to the idempotent $e$.
			\end{Lemma}

	\section{Algebraic Bernstein algebras}
	
	\label{sec:alg_Bernstein_algebras}
	
	Let $A$ be an arbitrary non-associative algebra over $K$. An element $a\in A$ is said to be {\it algebraic} if
	the subalgebra $\text{alg}(a)$ of $A$ generated by $a$ is finite-dimensional.  In this case, we set $\deg(a)=\dim \text{alg}(a)$, called {\it the degree} of $a$. The algebra $A$
	is said to be {\it algebraic} if every element of $A$ is algebraic. If, in addition, the algebraic degrees of all elements
	of $A$ are collectively bounded, we say that $A$ is {\it  algebraic of bounded degree}. In this case, the integer $\min\{n\in\N^*\;|\; \deg(a)\leq n\;\ \forall \,a\in A\}$ is called  {\it the (bounded) degree} of $A$ and is denoted by $\deg(A)$. Otherwise, we write
	$\deg(A)=\infty$.\\

	\esp Now, let $A$ be a Bernstein algebra. Since $A$ is not necessarily
	power-associative, the definition of the  powers of an
	element $a\in A$ has several variants. But here we use the principal powers  given recursively by
	$a^1=a,\dots, a^i=a^{i-1}a$. The following useful formula is well known for each
	element $a$ in a Bernstein algebra $A$ (see \cite{Koul},
	\cite{Wa1}):
	\begin{equation}
	\label{id:power_prod}
	a^ia^j=\frac12(\omega(a)^ia^j+\omega(a)^ja^i) \;\;\; \mbox{ for all } i, j\geq	2.
	\end{equation}
	This implies that an arbitrary power of $a$ (with any distribution
	of parentheses) is a linear combination of a finite number of its
	principal powers $a^i$. Hence, the subalgebra $\text{alg}(a)$ generated by $a$ is just the subspace of $A$ spanned
	by the principal powers of $a$, that is, $\text{alg}(a)=\lan a,
	a^2, a^3, \ldots\ran$. It follows that $a$ is algebraic if and only if there exist scalars $\gamma_1, \dots, \gamma_{n}\in K$ such
	that 
	\begin{equation} 
	\label{eq:a^n}
	a^{n+1}=\gamma_1 a+ \gamma_2 a^2 +\dots+ \gamma_{n}a^{n}
	\end{equation}
		 If $n$ is the least such number, then the elements $a, a^2, \dots, a^n$ are linearly independent, and so $\deg(a)=\dim \langle a, a^2, \dots, a^n\rangle =n$. It is noteworthy  that we do not include a constant term in the preceding relation, because a Bernstein algebra $A$ lacks a unity element. In other words, equation \eqref{eq:a^n} becomes $P(a)=0$ for some non-zero (ordinary) polynomial $P\in K[X]$ with $P(0)=0$.\\
		 On the other hand, for each element $x\in A$ and polynomial $Q\in K[X]$ with no constant term, it is not difficult to see that $\omega(Q(x))=Q(\omega(x))$.\\
		An element $a\in A$ is said to be {\it right nilpotent} if $a^m = 0$ for some $m\geq 2$, and the smallest such integer $m$ is called the {\it right nilpotency index} of $a$. Clearly, any right nilpotent element $a\in A$ of right nilpotency index $m$ is algebraic of degree $m-1$, since the subalgebra $\text{alg}(a)=<a, a^2, \dots, a^{m-1}>$ has dimension $m-1$. As remarked by Ouattara in \cite{Ouatt}, we will see in Section~\ref{sec:Jacobian} that a right nilpotent element of right nilpotency $m\geq 2$ in a Bernstein algebra is in fact {\it nilpotent} of the same {\it nilpotency index} $m$, in the sense that each product of $m$ copies of $a$, with any distribution of parentheses, vanishes.    
	
	\begin{Note}  \label{unit algebra} 
	{\rm  Let $(A, \omega)$ be an arbitrary Bernstein algebra. As evoked in \cite[Proposition 3.5.15]{Ly}, it is easily checked that for each $a\in H$ the subspace $\left(\text{alg}(a)\right)^2=<a^2, a^3, \dots >$ is an elementary Bernstein algebra, i.e. satisfying the identity $x^2=\omega(x)x$.}
	\end{Note}
	
	\begin{Note}\label{minimal polynomial}
	{\rm Let $a$ be algebraic as in (\ref{eq:a^n}). Classically, the set  $I_a$ of polynomials that have no constant term  and annihilate $a$ is an ideal in $K[X]$, and is therefore generated by a unique monic polynomial, denoted by $p_a$ and called {\it the minimal polynomial of} $a$. Thus, $p_a$ divides any other polynomial in $K[X]$ without a constant term that annihilates $a$.\\Indeed, take $p(X)=\sum\limits_{i=1}^n \alpha_iX^i\in I_a$ and let $q(X)=\sum\limits_{j=0}^s \beta_jX^j\in K[X]$.
Then $(pq)(X)=\sum\limits_{j=0}^s \beta_j\left(\sum\limits_{i=1}^n \alpha_iX^{i+j}\right)$ implies 
$$(pq)(a)=\sum\limits_{j=0}^s \beta_j\left(\sum\limits_{i=1}^n \alpha_ia^{i+j}\right)=\sum\limits_{j=0}^s \beta_j L_a^j\left(\sum\limits_{i=1}^n \alpha_ia^{i}\right)=\sum\limits_{j=0}^s \beta_jL_a^j(p(a)).$$ Thus, $(pq)(a)=0$ because  $p(a)=0$. \\
Notice that $\deg p_a=\deg(a)+1$ and $\text{alg}(a)=\{Q(a)\; /\; Q(X)\in K[X],\ \deg(Q)<\deg(p_a)\}$.
We will later establish in Proposition \ref{prop:min_pols_for_algebraic_elements} that if $\deg(a)=n\geq 2$, we have $\gamma_1=0$, so  (\ref{eq:a^n}) is reduced to
	\begin{equation} 
	\label{a^n}
	a^{n+1} = \gamma_2 a^2 +\dots+ \gamma_{n}a^{n},
	\end{equation}
and $p_a(X)=X^{n+1}-\gamma_nX^n -\dots -\gamma_2X^2$.
We also show in Proposition \ref{prop:min_pols_for_algebraic_elements} that the minimal polynomial $p_a$ of an algebraic element $a$ admits some explicit  form.}
	\end{Note}
		\esp Before beginning the investigation of algebraic Bernstein
	algebras, let us exhibit some concrete examples with regard to their algebraicity:
	
	\begin{Example} 
	\label{ex:fdim_Bernstein}
	{\rm Finite-dimensional Bernstein algebras are obviously algebraic of bounded degrees.}
	\end{Example}
	
	\begin{Example}
		\label{ex:Bernstein-Jordan}
		{\rm Bernstein Jordan algebras are algebraic of bounded degrees, since they satisfy the identity $x^3=\omega(x)x^2$.
		}
	\end{Example}
	
	\begin{Example} 
		\label{ex:train=>algebraic}
		\label{ex:Bernstein_train_is_algebraic}
		{\rm More generally, every Bernstein  train algebra of rank $r$ is algebraic of bounded degree $r-1$. Indeed, if the train equation of $A$ is $$x^r + \gamma_1\omega(x) x^{r-1} +
			\cdots+\gamma_{r-1}\omega(x)^{r- 1}x=0, \quad \mbox{  where } r=\text{rank}(A),$$ then for every element $a\in A$, we have $\dim \text{alg}(a)\leq r-1$, since $\text{alg}(a)$ is spanned by the elements $a, a^2, \ldots, a^{r-1}$. In this context, Lyubich proved in \cite[Corollary 3.5.13]{Ly} the existence of an element $a\in H$ such that $\dim \text{alg}(a)=r-1$.}
	\end{Example}
	
	Bernstein train algebras give a major source of algebraic Bernstein algebras. In the finite-dimensional case, we make use of the  crucial results of \cite{Ouatt}, \cite {Wa1} and \cite[Theorem 3.5.21]{Ly}, to which we add a new equivalent condition (iii). It will be later generalized (in Section~\ref{sec:train})  to the infinite-dimensional situation. More precisely:

\begin{Theorem} 
\label{th:2.5}
\label{th:f_dim_train}
Let $A=Ke\oplus U\oplus V$ be a finite-dimensional Bernstein algebra. The following
conditions are equivalent:\\
(i) $A$ is a train algebra;\\
(ii) For all $v\in V$, the operators $L_v: U \rightarrow U$ are  nilpotent. \\
(iii) For all $v\in V$, the operators $L_v: L(A)  \rightarrow L(A) $ are  nilpotent. \\
(iv) For all $x\in N$, the operators $L_x: N \rightarrow N$ are
nilpotent.\\
(v) $N$ is a nil-algebra.\\
In this case, the train equation of $A$ has the form
		$(x^3-\omega(x)x^2)(x-\frac 12\omega(x))^{r-3}=0$, where $r$ is the rank of $A$.
\end{Theorem}


	\begin{Note} {\rm Since our algebras need not be  power-associative, the above train equation $(x^3-\omega(x)x^2)(x-\frac 12\omega(x))^{r-3}=0$
			should be rigorously put in the form	
			$x^r+\gamma_1\omega(x)x^{r-1}+\dots+\gamma_{r-1}\omega(x)^{r-1}x=0$, where
			$$\gamma_k=\frac{(-1)^k}{2}
			\left[\binom{r-3}{k}+2\binom{r-3}{k-1}
			\right], \quad 1\leq k\leq r-1$$
			But this abusive notation is often used for the sake of simplification.}
			\end{Note}
	
	\begin{Example} 
	\label{ex:nuclear}
	{\rm  As a particular case of Bernstein train algebras, nuclear Bernstein algebras are
			algebraic of bounded degrees. Indeed, it is known from 		\cite[Th\'eor\`eme 3.1]{Mi1} (see also \cite[Proposition 9]{Catalan}) that such algebras are train algebras satisfying the train equation $x^4-\frac 32\omega(x)x^3+\frac 12\omega(x)^2x^2=0$. }
	\end{Example}
	
	The following example shows that the class of algebraic Bernstein algebras contains strictly the class of Bernstein-train algebras.
	\begin{Example}
	\label{ex:former_2.8}
	\label{ex:2.8}
	\label{ex:algebraic_not_train}
	{\rm There
			exist algebraic Bernstein algebras that are not train algebras. It is
			enough to consider  finite-dimensional Bernstein algebras. For instance, take the three-dimensional
			Bernstein algebra $A$ with basis $\{e, u, v\}$ and nonzero
			products $e^2=e, \; eu=\frac 12u, \; uv=u$. This algebra is not 
			train, since for $x=u+v\in N$, we have $x^3=x^2$, so $x$
			is not nilpotent and therefore $N$ is not nil. But the element $a=e+u+v\in A$ generates the whole algebra $A$, since the elements $a,\; a^2=e+3u,\; a^3 =e+5u$ are linearly independent. Thus, $A$ is algebraic of degree $3$.}
	\end{Example}
	
	\begin{Example} 
		\label{ex:former_2.9}
	\label{ex:Holgate_principal_train}
{\rm  Contrarily to Example~\ref{ex:train=>algebraic}, an arbitrary train algebra, which is not a
			Bernstein algebra, is not in general algebraic. To illustrate this
			situation, we make appeal to an example of Holgate taken from
			\cite{H}. Let $A$ be the {
			free train algebra
				corresponding to the train equation} $x^4-\omega(x)^3x=0$. Namely,
			$A={\cal F}_1/{\cal R}_{1,3}$, where ${\cal F}_1$ is the free
			commutative non-associative algebra without identity, with one
			generator $a$, endowed with the weight function such that
			$\omega(a)=1$, and ${\cal R}_{1,3}$ is the ideal of ${\cal F}_1$
			generated by all elements of the form $x^4-\omega(x)^3x, \;x$
			running over ${\cal F}_1$. Holgate proved that $A$ is
			infinite-dimensional, which shows that $a$ is not algebraic.\\
	However, power-associative train algebras are
	algebraic.  A systematic study of
	power-associative train algebras has been done in \cite{Ouatt1}. Especially, finitely generated Jordan train algebras were established to be finite-dimensional. }
\end{Example}
	
	\begin{Example} 
	\label{ex:non-algebraic}
	{\rm There are several Bernstein algebras that are not algebraic.
			For instance, let us consider the Bernstein algebra $A=Ke\oplus
			U\oplus V$ treated in \cite[Example 3.10]{B} and  \cite{Su},
			where $U=\langle u_1, u_2, \ldots\ran ,\; V=\langle v\ran$ and the multiplication table
			given by
			$$e^2=e,\; eu_i=\frac 12 u_i, \;
			u_iv=u_{i+1} \; (i\geq 1),\;\; \mbox{  other products being zero. }$$
			Then, $A$ is not algebraic. Indeed, the principal powers of the element
			$x=u_1+v$ are $x^2=2u_{2}, \,
			x^3=2u_{3}, \dots, x^{i}=2u_{i}$  for all $i\geq
			2$. It follows that $\text{alg}(x)=\lan u_1+v, u_2, u_3, \ldots\ran $ is infinite-dimensional, so $x$ is not algebraic.}
	\end{Example}

	\begin{Example}
	\label{ex:former_2.11}
	\label{ex:algebraic_unbounded_degree}
	{\rm  An algebraic Bernstein algebra is not necessarily of bounded
			degree. Indeed, let $A=Ke\oplus U\oplus V$ be the Bernstein
			algebra considered in \cite[Example 3.10]{B}, \cite[Example
			4.6]{Mi},  where $U=\langle u_1, u_2, \ldots\rangle , \;V=\langle v\ran$ and the
			multiplication table is given by
			$$e^2=e,\; eu_i=\frac 12 u_i \; (i\geq 1),  \;
			u_iv=u_{i-1} \; (i\geq 2),\;\; \mbox{ other products being zero.}$$
			Let $x=\alpha e+\sum\limits_{i=1}^n \beta_i u_i+\lambda v\in A$.
			Then an easy   induction shows  that all powers $x^n$ belong to the finite-dimensional
			subspace $T_n=\langle e, v, u_1, u_2, \ldots, u_n\rangle $. So, $\text{alg}(x)\subseteq
			T_n$, and therefore $A$ is algebraic. Nevertheless, $A$ is not of
			bounded degree. For this, taking the element
			$x_n=u_n+v$ with $n\geq 1$, we have $x_n^2=2u_{n-1}, \,
			x_n^3=2u_{n-2}, \dots, x_n^{n}=2u_1$ and $x_n^{i}=0$ for $i\geq
			n+1$. It follows that $\text{alg}(x_n)=\langle u_1, u_2, \dots, u_{n-1}, u_n+v\rangle \,$
			and $\dim \text{alg}(x_n)=n$. Therefore, $A$ is algebraic with unbounded
			degree. \\
			This example also shows that
			infinite-dimensional algebraic Bernstein algebras need not be
			train algebras, see Remark~\ref{ex:locally_train_non_train} below.
			}
	\end{Example}

	\esp We now recall some definitions on linear operators. Let $E$ be an arbitrary vector space over $K$ and $T: E
	\rightarrow E$ be a linear operator. Recall that $T$ is {\it
		algebraic} if there is a non-zero (ordinary) polynomial $P(X)\in
	K[X]$ such that $P(T)=0$. In this case, we write
	$\deg(T):=\deg(P)$ and call it the {degree} of $T$, where $P(X)$
	is the polynomial of least degree annihilating $T$. The operator
	$T$ is said to be {\it locally algebraic} if, for every element
	$a\in E$, there is a non-zero polynomial $P_a(X)\in K[X]$ such that
	$P_a(T)(a)=0$. If, in addition, the degree $\deg(P_a)$ is bounded
	independently of $a$, we say that $T$ is locally algebraic {\it
		with  bounded degree}.\\

	\esp Our first result will formulate the algebraicity of a
	Bernstein algebra $A$ in terms of some multiplication operators.
	Before  starting discussion, we require some preparation by proving
	the next auxiliary lemma.
	
	\begin{Lemma} 
	\label{lem:Lemma_2.12}
	\label{lem:former_Lemma_2.9}
	    Let $A=Ke\oplus U\oplus V$ be a Bernstein algebra.
		For all $u\in U, v\in V$ and $k\geq 2$, we have
		$$  (u^2+v)^k=2L_v^{k-1}(u^2)+v^k, $$
		where $L_v$ is the multiplication operator by $v$.
	\end{Lemma}
	
	{\it Proof. \ } It is clear for $k=2$, since
	$(u^2+v)^2=(u^2)^2+2u^2v+v^2=2L_v(u^2)+v^2$.\\
	Let $k\geq 2$. By induction,
	$(u^2+v)^{k+1}=(u^2+v)(2L_v^{k-1}(u^2)+v^k)=2u^2L_v^{k-1}(u^2)+u^2v^k+2L_v^{k}(u^2)+v^{k+1}$.
	In view of (1.5), we have $u^2v^k=0$ for $k\geq 2$ (see, also, \cite[(1.14)]{Ouatt} or
	\cite[Lemma 3.1,(1)]{GarciaGonzalez}). Moreover, $VU^2\subseteq VV=V^2\subseteq L(A) $, and since $L(A) $ is an ideal, we get $V(V(\dots V(VU^2)\dots))\subseteq L(A) $. Hence, $L_v^{k-1}(u^2)\in L(A) $, so that $u^2L_v^{k-1}(u^2)=0$, because $U^2.L(A) =0$ by (1.7).
	Therefore, $(u^2+v)^{k+1}=2L_v^k(u^2)+v^{k+1}$.  \ep\\
	
	\esp Now, we are ready to establish our first principal result.
	
	\begin{Theorem} 
	\label{th:A_is_alg_iff_Lv_is_alg}
	\label{th:former_Th2.13}
	Let $A=Ke\oplus U\oplus V$ be a Bernstein
		algebra and let $N=\ker(\omega)=U\oplus V$. Then the following
		conditions are equivalent:\\
		(i) The algebra $A$ is algebraic;\\
		(ii) The ideal $N$ is algebraic;\\
		(iii) The operators $L_v: U \rightarrow U$ are locally
		algebraic for all $v\in V$.
	\end{Theorem}
	
	{\it Proof. \ }   (iii) $\Rightarrow$ (i): Let $x=\alpha e+u+v\in A$ and $k\geq 4$. By
	\cite[Proposition 3.4]{GarciaGonzalez},
	$$x^k=\alpha^ke+ \left(\alpha^{k-1} u+ b_1^{(k)}\alpha^{k-2}L_v(u)+\dots
	+b_{k-1}^{(k)}L_v^{k-1}(u)\right)+\left(a_2^{(k)}\alpha^{k-2}
	v^2+\dots+ a_k^{(k)} v^k\right)$$
	$$+\left(\alpha^{k-2}u^2+c_3^{(k)}\alpha^{k-3}L_v(u^2)+c_2^{(k)}\alpha^{k-4}L^2_v(u^2)+\dots+c_{k-2}^{(k)}L^{k-2}_v(u^2)\right)$$
	for some scalars $a_i^{(k)}, b_i^{(k)}$ and $c_i^{(k)}$ in $K$.
	Therefore, $x^k$ belongs to the subspace $S_{u, v}$ given by
	$$S_{u, v}=\langle e,\; L_v^i(u)\ (i\geq 0), \; L_v^j(v)\ (j\geq 1), \; L_v^k(u^2)\ (k\geq 0)\rangle  
	$$
	$$\qquad=
	\;\langle e,\; L_v^i(u)\ (i\geq 0), \; L_v^i(v^2)\ (i\geq 0), \;
	L_v^i(u^2)\ (i\geq 0)\rangle ,
	$$
	where $v^2\in U$ because $V^2\subset U$.
	Now, since the operator $L_v: U
	\rightarrow U$ is locally algebraic, the subspaces $\langle L_v^i(u) \
	(i\geq 0)\rangle , \; \langle L_v^i(v^2) \ (i\geq 0)\rangle  $ and $\langle L_v^i(u^2) \
	(i\geq 0)\rangle = \langle u^2, \; L_v^i(vu^2) \ (i\geq 1)\rangle $ of $A$ are
	finite-dimensional, and so is $S_{u, v}$.
	\\Finally, since
	$\text{alg}(x)=\langle x, x^2, \dots, \rangle $, we deduce that $x$ is algebraic.\\
	Observe that $x^2=\alpha^2+\alpha u+2uv+u^2+v^2$ and
	$x^3=\alpha^3e+2\alpha(uv)+\frac 12\alpha v^2+\alpha^2
	u+\alpha^2u^2+2u(uv)+uv^2+u^3+2(uv)v+v^3+u^2v=\alpha^3e+\alpha^2u+2\alpha(uv)+\frac
	12\alpha v^2+2(uv)v+v^3$ (see the proof of \cite[Lemme
	2.5]{Ouatt}). Thus, $x^k\in S_{u, v}$  for any $k\geq 1$ and not only for $k\geq 4$. \\
	
	(ii) $\Rightarrow$ (iii): Let $u\in U$ and $v\in V$.
	Since $u^2+v$ and $v$ are both algebraic of degrees $\leq n$, the
	subspaces $F_1=\langle (u^2+v)^k \; |\; k\geq 1\rangle $ and $F_2=\langle v^k \; |\;
	k\geq 1\rangle $ are finite-dimensional. Now, by Lemma~\ref{lem:former_Lemma_2.9}, we have
	$L_v^{k-1}(u^2)\in F_1+F_2$ for all $k\geq 1$. So, the subspace
	$\langle L_v^{k-1}(u^2) \; |\; k\geq 1\rangle $ is finite-dimensional, too.\\
	On the other hand, in virtue of \cite[Proposition 3.4]{GarciaGonzalez}, we
	have, for $x=u+v$ and $k\geq 4$,
	\begin{equation} 
	\label{eq:x^k}
	x^k=a_k^{(k)} v^k+b_{k-1}^{(k)}
		L_v^{k-1}(u)+c_{k-2}^{(k)} L_v^{k-2}(u^2),
	\end{equation}
	for some coefficients $a_k^{(k)}, \, b_{k-1}^{(k)}, \,
	c_{k-2}^{(k)}$ in $K$, with $b_{k-1}^{(k)}=2$ \cite[Lemma 3.1,
	(1)]{GarciaGonzalez}. Since $x$ and $v$ are algebraic, the subspaces
	$\langle x^k\;|\; k\geq 1\rangle $ and $\langle v^k\; |\; k\geq 1\rangle $ are
	finite-dimensional. In addition, we have already seen that the
	subspace $\langle L_v^{k-1}(u^2)\; |\; k\geq 1\rangle $ is also finite-dimensional.
	We conclude from (\ref{eq:x^k}) that the subspace $\langle L_v^{k-1}(u)\;|\; k\geq
	1\rangle $ is finite-dimensional, for all $u\in U$ and $v\in V$. As a
	consequence, the operator $L_v$ is locally algebraic for every $v\in 	V$.
	
	Since the implication (iii) $\Rightarrow$ (i) is obvious, the proof is complete. 
	\ep\\
	
	 The next result by Kaplansky  (see \cite{Aupetit}) will be useful in our
	development. 
	\begin{Lemma} {\rm(Kaplansky)} Let $E$ be a vector space and $T: E \rightarrow
		E$ be a linear operator. If $T$ is locally algebraic of bounded
		degree, then $T$ is algebraic.
	\end{Lemma}
	\esp With the help of this lemma, we may study whether a Bernstein algebra is algebraic of bounded degree, completing Theorem~\ref{th:A_is_alg_iff_Lv_is_alg}.
	
	\begin{Theorem} 
	\label{th:Lv_is_alg_bounded}
	\label{th:Theorem_2.15}
		 For a Bernstein algebra $A=Ke\oplus U\oplus V$ with  barideal  $N=\ker(\omega)=U\oplus V$,
	   the following
		conditions are equivalent:\\
		(i) The algebra $A$ is algebraic of bounded degree;\\
		(ii) The ideal $N$ is algebraic of bounded degree;\\
		(iii) The operators $L_v: U \to U$ are algebraic for all $v\in V$ and their
		algebraic degrees are bounded collectively.
	\end{Theorem}
	
	{\it Proof. \ } Let us return to the proof of Theorem~\ref{th:A_is_alg_iff_Lv_is_alg} and use the
	same notation. It is sufficient to prove the implications (ii) $\Rightarrow $ (iii) $\Rightarrow $ (i).
	
	(ii) $\Rightarrow $ (iii): Assume $A$ is algebraic of
	bounded degree $n$. By the implication (i) $\Rightarrow$ (ii) of Theorem~\ref{th:A_is_alg_iff_Lv_is_alg} and
	its proof, the operators $L_v$ are locally algebraic for all $v\in
	V$. More specifically, $\dim F_1\leq n$ and $\dim F_2\leq n$, so
	that $\dim \langle L_v^{k-1}(u^2) \; |\; k\geq 1\rangle  \leq 2n$. It follows
	that $\dim \langle L_v^{k-1}(u) \; |\; k\geq 1\rangle  \leq 4n$ for all $u\in U$
	and $v\in V$, which implies that $L_v$ is locally algebraic of
	bounded degree $\leq 4n$, for all $v\in V$. Finally, according to
	Kaplansky's lemma, we infer that $L_v$ is algebraic of degree
	$\leq 4n$, for all $v\in V$.\\
	
	(iii) $\Rightarrow$ (i): By hypothesis, the algebraic degrees of
	the $L_v$ are bounded independently of $v$, that is, there exists
	an integer $m$ such that $\dim \text{alg}(L_v)\leq m$ for all $v\in V$.
	Consequently, the dimension of the subspace $S_{u, v}$ is bounded
	independently of $u$ and $v$, and so is the subspace $\text{alg}(x)$
	(which satisfies $\dim \text{alg}(x)\leq m+2$). \ep\\
	
	\esp  The next example shows that it is indispensable that $A$ has bounded degree in Theorem~\ref{th:Lv_is_alg_bounded}.
	
	\begin{Example} {\rm Let us consider again our previous Example~\ref{ex:algebraic_unbounded_degree}. Since $L_{v}^k(u_i)=u_{i-k}$ for all $i\geq 2$ and $k\geq 1$ with $i>k$, then the operator $L_{v}$ cannot be annihilated by a polynomial. Thus, $L_{v}$ is not algebraic (but only locally algebraic), although the algebra $A$ is algebraic.  }
	\end{Example}

	\esp Actually, the study of Bernstein algebras that are
	algebraic is equivalent to the study of the associated linear
	operators $L_v$, and this approach will be important in  Section~\ref{sec:Kurosh} 
	about
	the Kurosh problem.\\
	
\esp We now continue to establish some  properties of algebraic Bernstein algebras.
	Our next task is to reduce the problem to the Lyubich ideal $L(A) $.
	
 Recall that a Bernstein algebra $A$
	is Jordan if $A$
	satisfies the identity  $x^3-\omega(x)x^2=0$.
	Hence, if $I$ is a baric ideal of $A$ such that the quotient
	algebra $A/I$ is Jordan, then $x^3-\omega(x)x^2\in I$ for all
	$x\in A$. In particular, $x^3-\omega(x)x^2\in L(A) $ for all
	$x\in A$. Let us denote by $l_v: L(A)  \rightarrow L(A) $
	the restriction of $L_v$ to $L(A) $. The following result
	refines Theorems~\ref{th:A_is_alg_iff_Lv_is_alg} and~\ref{th:Lv_is_alg_bounded}.
	
	\begin{Proposition} 
	\label{pr:2.18}
	For a Bernstein algebra $A=Ke\oplus U\oplus V$, the following
		conditions are equivalent:\\
		(i) The algebra $A$ is algebraic (of bounded degree);\\
		(ii) The subalgebra $B_e=Ke\oplus L(A) \oplus V$ is algebraic
		(of
		bounded degree);\\
		(iii) The operators $l_v: L(A) 
		\rightarrow L(A) $ are locally algebraic (algebraic) for all $v\in V$
		(and their algebraic degrees are bounded collectively).
	\end{Proposition}
	
	{\it Proof. } The implication (i) $\Rightarrow$ (ii) is obvious
	and the equivalence (ii) $\Longleftrightarrow$ (iii) is immediate
	by applying Theorems~\ref{th:A_is_alg_iff_Lv_is_alg} and~\ref{th:Lv_is_alg_bounded} to the Bernstein algebra $B_e$. \\
	(iii) $\Rightarrow$ (i): Let $x=\alpha e+u+v\in A$. By \cite[(2.6)]{Ouatt}, we have
	\begin{equation}
	\label{eq:LHS_train_pol_via_Lv}
	\left( x^3-\omega(x)x^2  \right) \left(x-\frac 12\omega(x)\right)^k=L_v^k(x^3-\omega(x)x^2), \quad \mbox{ for all }k\geq 1.
	\end{equation}
	Since $l_v$ is locally algebraic and $x^3-\omega(x)x^2\in L(A) $,
	there exist scalars $a_0, \dots, a_p\in K$, not all zero, such that
	$\displaystyle{\sum\limits_{k=0}^p a_k l_v^k(x^3-\omega(x)x^2)=0}$, which
	gives $$\displaystyle{\sum\limits_{k=0}^p a_k (x^3-\omega(x)x^2)\left(x-\frac
		12\omega(x)\right)^k =0}.$$ Hence, $x$ is algebraic and the proof is achieved.
	\ep
	\begin{Note}{\rm We may refine condition (ii) in the previous proposition by furnishing a weaker condition. Let $\overline{A}=Ke\oplus L(A) \oplus V$ be the new algebra with the same multiplication as subalgebra of $A$ with the exception of the products $vv'$ that we assign to be zero for all $v, v'\in V$. Clearly, $\overline{A}$ is a Bernstein algebra by Lemma~\ref{lemma:Lemma1.1}. It is immediate from the equivalence (ii) $\iff$ (iii) that  the Bernstein algebra  $A$ is  algebraic (of bounded degree) if and only if the Bernstein algebra $\overline{A}$ is. This means that the subspace $V^2$ in a Bernstein algebra has no influence on the fact that $A$ be algebraic or not.}	
	\end{Note}
	
	\begin{Proposition} \label{alg H} Let  $A=Ke\oplus U\oplus V$ a Bernstein algebra with $H = \{x\in A \;|\; \omega(x)=1 \}$. The following
		conditions are equivalent:\\
		(i) The algebra $A$ is algebraic (of bounded degree);\\
		(ii) All elements $x\in H$
		are algebraic (and their algebraic degrees are bounded collectively).
	\end{Proposition}
	
	\pr 
	The implication (i)$\Rightarrow$(ii) is obvious in both the bounded and unbounded degree cases.  
	Let us prove (ii)$\Rightarrow$(i).
	
	Assume that each element $x\in H$ is algebraic. It follows that all $x\in A$ with $\omega(x)\neq 0$ are algebraic.  Let us prove that $A$ is algebraic. 
	 Using Theorem~\ref{th:A_is_alg_iff_Lv_is_alg}, we show that the operators $l_v: L(A) 
	\rightarrow L(A) $ are then locally algebraic for all $v\in V$.\\
	
	{\it First step}: For $x=e+v$ and $k\geq 4$, we have by \cite[Proposition 3.4]{GarciaGonzalez}:
	$$x^k=e+\left(a_2^{(k)}\alpha^{k-2}
	v^2+\dots+ a_k^{(k)} v^k\right) \mbox{ with } a_k^{(k)}=1.$$
	It follows that $v^2\in\ \langle e, x^2\rangle ,\, v^3\in\ \langle e, v^2, x^3\rangle  \subset \langle e, x^2, x^3\rangle $, and by induction \newline $v^k\in\ \langle e,x^2, x^3, \dots, x^k\rangle $. Since $x\in H$ is algebraic,  the subspace $\langle e,x^2, x^3, \dots, x^k\rangle $ is finite-dimensional, and so also is the subspace $\langle v, v^2, v^3, \dots\rangle $.\\
	
	{\it Second step}: Let $y=e+u+v\in A$ with $u\in L(A) $ and $v\in V$. Applying again \cite[Proposition 3.4]{GarciaGonzalez}, and using $u^2=0$, we get
	$$y^k=e+ \left(u+ b_1^{(k)}L_v(u)+\dots
	+b_{k-1}^{(k)}L_v^{k-1}(u)\right)+\left(a_2^{(k)}
	v^2+\dots+ a_k^{(k)} v^k\right)$$
	Thus, it easy to prove by induction that for all $i\geq 1,$
	$$L_v^i(u)\in \ \langle e, u\rangle + \langle y^k\, |\, k\geq 2\rangle +\langle v^k\, |\, k\geq 2\rangle $$
	It follows from the preceding arguments that the subspace $\langle L_v^i(u)\,|\, i\geq 1\rangle $ is finite-dimensional, since the subspaces $\langle y^k\, |\, k\geq 1\rangle $ and $\langle v^k\,|\, k\geq 1\rangle $ are finite dimensional. \\
	In conclusion, the operator $L_v: L(A)  \rightarrow L(A) $ is locally algebraic for all $v\in V$. 
	This proves the unbounded version of the implication (ii)$\Rightarrow$(i).\\
	
	Assume now that all elements  $x\in H$ are algebraic and their algebraic degrees are bounded collectively by some $D>0$. It follows that  all $x\in A$ with $\omega(x)\neq 0$ are algebraic, and their 
	algebraic degrees are bounded by $D$ as well.  Let us prove that each element $n\in N$ is algebraic of degree at most $D$. For each $x_\alpha = \alpha e+n$ with $\alpha \in K^\times$, the powers 
	$\{ x_\alpha^k | k =1,\dots, D\} $ are linearly dependent, where $x_\alpha^k = n^k + \alpha (y_{\alpha,k})$
	for some $y_{\alpha,k}\in A$.
	It follows that there exist polynomials $P_1(\alpha), \dots, P_D(\alpha)$ (some of them are non-zero) 
	such that  $\sum_{k=1}^{D} P_k(\alpha) x_\alpha^k = 0$.
	One may assume that for some $j\in [1,D]$ we have $P_j(0)\ne 0$. Then we get 
 a linear dependence $\sum_{k=1}^{D} P_k(0) x_0^k = 0$ of the elements $n^k = x_0^k$.
	  \ep\\

	\begin{Theorem} 
	\label{th:alg_ideal}
	Any Bernstein algebra $A$ contains a maximal algebraic ideal ${\cal I}(A)$ of $A$. Furthermore, ${\cal I}(A)$ contains $A^2=Ke\oplus U\oplus U^2$
	and $A/{\cal I}(A)$ is a zero-multiplication algebra. 
	\end{Theorem}	
	
	\pr \quad Following Levitzki \cite[page 404]{Levitzki}, if $B$ and $C$ are algebraic ideals in a Bernstein algebra $A$, let us prove that the ideal $B+C$ is also algebraic. Let $b\in B$ and $c\in C$. An easy induction shows that for each $k\geq 1$, there exists $c_k\in C$ such that $(b+c)^k=b^k+c_k$.  Now, let $P(X)=\sum\limits_{k=1}^m \alpha_k X^k \in K[X]$ with $P(b)=0$ and $\deg(P)\geq 1$. Then $P(b+c)=f(b)+\sum\limits_{k=1}^m \alpha_k c_k=c'$, where $c'\in C$. Let $Q(X)=\sum\limits_{k=1}^n \beta_k X^k \in K[X]$ with $Q(c')=0$ and $\deg(Q)\geq 1$. Thus, $Q(P(b+c))=0$. Observe that for a lack of power-associativity, the equality $Q(P(b+c))=(Q\circ P)(b+c)$ need not be true. Nonetheless, we may prove that $Q(P(b+c))=R(b+c)$ for some polynomial $R\in K[X]$. Indeed, it is enough to use the identity~(\ref{id:power_prod}) in the form 
	\begin{equation}
	\label{id:power_prod_xy}
	  x^i x^j = \frac 12\left[y^i x^j+y^jx^i\right],
	\end{equation}
	where $x = b+c$, $y = \omega(b+c)$, and $i,j \ge 2$.
	
	It remains to check that the polynomial $R$ is nonzero. This follows from Lemma~\ref{lem:xy_powers} below.
	
	On the other hand, the subalgebra $A^2$ is an ideal of $A$.  
	Further, according to Example~\ref{ex:nuclear}, $A^2$ is train algebra. Consequently, $A^2$ is algebraic, which yields $A^2\subset {\cal I}(A)$.\\
	We may write the ideal ${\cal I}(A)$ in the form ${\cal I}(A)=Ke\oplus U\oplus V'$, where $V'$ is a subspace of $V$ with $U^2\subset V'$. This comes also from \cite{Koul}, since the ideal  ${\cal I}(A) \not\subset N$ ($e\in {\cal I}(A)$). Then the multiplication in the quotient algebra $A/{\cal I}(A) \simeq V/V'$ is trivial, since $V^2\subset U \subset {\cal I}(A)$.   	\ep
	
	\begin{Lemma}
	\label{lem:xy_powers}
	Let $\mathcal{R}$ be the ring of non-associative commutative polynomials in the variable $x$ over the commutative associative ring $K[y]$ subject 
	to the relations~(\ref{id:power_prod_xy}) for all $i,j\ge 2$. Assume that $P$ and $Q$ are two polynomials in $x$  with coefficients in $K$ without constant terms. Then the composition polynomial $Q(P(x))$ is non-zero in $\mathcal{R}$.
	\end{Lemma}
	
	\pr
	The ring $\mathcal{R}$ is $\Z$-graded with the standard grading such that $\deg (x^i y^j) = i+j$. If the leading terms of the polynomials $Q$ and $P$ are $\alpha_m x^m $ and $\beta_n x^n$ respectively, then the only term of degree $mn$ in $R(x) = Q(P(x))$ is the term $M = \alpha_m \beta_n^m (x^n)^m$. Let $\phi: \mathcal{R}\to K[y]$ be the epimorphism which sends simultaneously $x$ and $y$ to $y$. Then $\phi(M) = \alpha_m \beta_n^m y^{mn} \ne 0$, so that 
	$M\ne 0$ and $R(x)\ne 0$.
	\ep
	
	
	\begin{Note}
	\label{rem:former_2.22}
	{\rm The inclusion $A^2\subset {\cal I}(A)$ is 
	in general strict, for example, for each finite-dimensional Bernstein algebra $A$ such that $U^2 \ne V$ (such as  finite-dimensional one-generated Bernstein algebras described in Theorem~\ref{th:principal_subalgebra} (ii) below). 
	
	In contrast to this situation, if $A$ is the infinite-dimensional one-generated algebra given below in Theorem~\ref{th:principal_subalgebra} (iii), then $A$ is not algebraic, $A^2 = Ke \oplus U$ and $V$ 
	is one-dimensional, so that  $ {\cal I}(A)=A^2$.} 
	\end{Note}	
	

	\section{Singly generated Bernstein algebras}
	
		\label{sec:singly_generated}

	In connection with the concept of algebraic elements, it would be desirable to investigate the structure of  singly generated subalgebras in  an arbitrary Bernstein algebra $A$. In this section we deal with this question by analyzing the Peirce decomposition of such subalgebras. As a consequence, we will be able to provide some information about the minimal polynomial $p_a$ of an algebraic element $a\in A$.\\  
	Obviously, the subalgebra $\text{alg}(x)$ for $x\in A$ is Bernstein if and only if $\omega(x)\neq 0$. Lyubich \cite[Corollary 3.5.18]{Ly} proved that in a finite-dimensional Bernstein algebra $A$, 
	if $x\in H$  is not an idempotent, then the Bernstein subalgebra $\text{alg}(x)$   has type $(n-1,1)$ and is exceptional. Subsequently, using a different approach, Conseibo and Ouattara  recovered this result \cite[Lemme 2.4]{Cons}  and furnished the multiplication table with respect to an appropriate basis of $\text{alg}(x)$ \cite[Théorème 3.9]{Cons}. 
	We revisit below these results and complete them for the infinite-dimensional context, by providing a detailed proof  for the sake of completeness.\\
	
If $x\in A$ with $\omega(x)\neq 0$, there is no loss of generality to assume $\omega(x)=1$, that is, $x\in H$, so that $e=x^2$ is a nonzero idempotent. First, let us look for Bernstein algebras of dimension $n\leq 3$ which are singly generated. Evidently, if $n=1$, then $x$ is itself an idempotent and $\text{alg}(x)=Ke$. Now, for $n=2$ or $3$, the following result was demonstrated in \cite[Propositions 2.3 and 3.6]{Cons}:

\begin {Proposition} \label{Conseibo} (i) The unique singly generated Bernstein algebra of dimension $2$ is  {\it the constant algebra} of type $(1,1): A_1=\langle e, v\rangle$ with $e^2=2, ev=v^2=0$.\\
(ii) Let $A$ be a Bernstein algebra of dimension $3$ generated by an element $x$ with $\omega(x)=1$.  Then the idempotent $e=x^2$ induces a Peirce decomposition $A=Ke\oplus U_x\oplus V_x$, where $U_x=\langle u_1\rangle$ and $V_x=\langle v_1\rangle$ with $u_1=x^3-x^2$ and $v_1=x+x^2-2x^3$. Moreover, $x=e+2u_1+v_1$ and there exists $\alpha\in K$ such that the multiplication table  is:
$$e^2=e, \; eu_1=\frac 12u_1,\; v_1^2=4(1-\alpha)u_1, \; u_1v_1=\Bigl(\alpha-\frac 32\Bigr)u_1, \; u_1^2=0.$$
\end{Proposition}


The next result concerns singly generated Bernstein algebras of dimensions $n\geq 4$.

\begin{Theorem} 
\label{th:principal_subalgebra}
Let $A$ be an arbitrary Bernstein algebra, and let $x\in H$. Consider  the principal Bernstein subalgebra 
$\text{\rm alg}(x)$.\\
(i) The idempotent $e=x^2$ induces a Peirce decomposition $\text{\rm alg}(x)= Ke\oplus U_x\oplus V_x$, where $U_x=\langle x^3-x^2, x^4-x^2, \dots \rangle, \; U^2_x=0, \;V_x=\langle v_1\rangle $ with $v_1=x+x^2-2x^3$.\\
(ii) If $x$ is algebraic of degree $n\geq 4$, then 
 there is a finite basis $(u_1, \ldots, u_{n-2})$ of $U_x$  and scalars $\beta_1, \ldots, \beta_{n-2}$ such that  $u_1=x^3-x^2$ and $$v_1u_i=u_{i+1}\, (1\leq i\leq n-3),\ v_1u_{n-2}=\beta_1u_1+\cdots+\beta_{n-2}u_{n-2},\, v_1^2=-2u_1-4u_2.$$
In this case,  $\text{\rm alg}(x)$ is a train algebra if and only if $\beta_1= \dots=\beta_{n-2}=0$.  \\
(iii) If $x$ is not algebraic, then 
there is a denumerable basis $(u_1, u_2, \ldots)$ of $U_x$ such that
$$u_1=x^3-x^2,\;v_1u_i=u_{i+1} \,(i\geq 1),\, v_1^2=-2u_1-4u_2.$$
In this case, the subalgebra $\text{\rm alg}(x)$ is N\oe therian but 
not
Artinian.
\end{Theorem}

\pr  The proof of (i) will be implicitly contained in those of (ii) and (iii).\\
 (ii) Assume that $x$ is algebraic of degree $n$, so that alg$(x)=\langle x, x^2, \ldots, x^n\rangle $. Let alg$(x)=Ke\oplus U_x\oplus V_x$ be the Peirce decomposition of $\text{alg}(x)$ with respect to the idempotent $e=x^2$. Select $y=\sum\limits_{i=1}^n\alpha_ix^i\in \text{alg}(x)$ with $\alpha_i\in K$. By identity \eqref{id:power_prod}, we have
$$ey=\left[\alpha_2+\frac 12(\alpha_3+\cdots+\alpha_n)\right]x^2+\left(\alpha_1+\frac 12\alpha_3\right)x^3+\left(\frac{\alpha_4}{2}\right)x^4+\cdots+
\left(\frac{\alpha_n}{2}\right)x^n$$
Hence, $ey=0 \iff \alpha_4=\dots=\alpha_n=0$ and $ \alpha_1=\alpha_2=-\frac 12\alpha_3 \iff y=\alpha_1(x+x^2-2x^3)$.
Thus, $V_x=\langle v_1\rangle $ is one-dimensional, with $v_1=x+x^2-2x^3$. Now, $ey=\frac 12y\iff \alpha_1=0$ and $\alpha_2=-(\alpha_3+\dots+\alpha_n)\iff
y=\alpha_3(x^3-x^2)+ \dots+\alpha_n(x^n-x^2)$. Therefore, $U_x=\langle x^3-x^2, x^4-x^2, \dots, x^n-x^2\rangle $ has dimension $n-2$. Observe that $x=e+u_x+v_x$, where $u_x=2(x^3-x^2)\in U_x$ and $v_x=v_1\in V_x$. Now, we have $U_x^2=0$, because for $1\leq i, j\leq n-2$,  we get by \eqref{id:power_prod}:
$$(x^{i+2}-x^2)(x^{j+2}-x^2)=\frac 12(x^{i+2}+x^{j+2})-\frac 12(x^{i+2}+x^2)-\frac 12(x^2+x^{j+2})+x^2=0.$$
On the other hand, we also have $U_x=\langle x^3-x^2, x^4-x^3, \dots, x^n-x^{n-1}\rangle$. Putting $u_1=x^3-x^2$ and $u_{i+1}=u_iv_1$ for $1\leq i\leq n-3$, we proceed by induction to show that 
\begin{equation}
	\label{eq:***} \tag{***}
x^{i+2}-x^{i+1}\in \langle u_1, \ldots, u_i\rangle, \;\;\;\mbox{ for } 1\leq i\leq n-2
	\end{equation}
This is clear for $i=1$, since $x^3-x^2=u_1$. Let $i\leq n-3$ and assume that $x^{i+2}-x^{i+1}=\lambda_1u_1+\dots+\lambda_{i}u_{i}$ with $\lambda_1, \dots, \lambda_{i}$ in $K$. Then 
$$x^{i+3}-x^{i+2}=x(x^{i+2}-x^{i+1})=(e+u_x+v_1)(\lambda_1u_1+\dots+\lambda_{i}u_{i})$$
$$=\frac 12(\lambda_1u_1+\dots+\lambda_{i}u_{i})+\lambda_1(v_1u_1)+\dots+\lambda_{i}(v_1u_{i})=\frac 12(\lambda_1u_1+\dots+\lambda_{i}u_{i})+\lambda_1 u_2+\dots+\lambda_{i} u_{i+1}$$
$\in \langle u_1, \dots, u_i, u_{i+1}\rangle$, as desired. Now, since 
$(x^3-x^2, x^4-x^3, \dots, x^n-x^{n-1})$ is a basis of $U_x$, it follows from  \eqref{eq:***} that  $(u_1, \dots, u_{n-2})$ is a spanning set for $U_x$, which is in turn a basis, because $\dim U_x=n-2$.\\
In addition, as $u_{n-2}v_1\in U_x$, there exist scalars $\beta_1, \ldots, \beta_{n-2}$ such that $u_{n-2}v_1=\beta_1u_1+\cdots+\beta_{n-2}u_{n-2}$. Finally, it is readily checked from  \eqref{id:power_prod} that $v_1^2=(x+x^2-2x^3)^2=4x^3-4x^4=-4x(x^3-x^2)=-4(e+u_x+v_1)u_1=-2u_1-4u_2$.\\
According to Theorem~\ref{th:2.5},  $\text{alg}(x)$ is a train algebra 
when
the operator
$L_{v_1}: U_x \rightarrow U_x$ is nilpotent.
But since this operator is cyclic and $L^{n-2}_{v_1}(u_1)=v_1u_{n-2}= \beta_1u_1+\cdots+\beta_{n-2}u_{n-2}$, its nilpotency yields 
$\beta_1= \cdots=\beta_{n-2}=0$.\\
(iii) If $x$ is not algebraic, then $\text{alg}(x)=\langle x, x^2, x^3, \dots\rangle $ is countably dimensional, and by the same above calculus we obtain the required conclusion.\\
Moreover, since $\text{alg}(x)$ is a finitely generated Bernstein algebra with $\dim V_x=1$, one derives from \cite[Proposition 3.7]{B} that $\text{alg}(x)$ is Noetherian. However, as the chain of ideals $I_k=\langle u_k, u_{k+1}, \dots\rangle  (k\geq 1)$ is strictly descending, $\text{alg}(x)$ is not Artinian.
\ep


\begin{Note}
\label{rem:single-generated-free}{\rm 
In the category of Bernstein algebras, the infinite-dimensional Bernstein algebra $A=\text{alg}(a)$ described in Theorem~\ref{th:principal_subalgebra}(iii) may be viewed as 
the  {\it free Bernstein algebra with one generator}, that is the largest Bernstein algebra generated by a single element. In effect, this algebra $A$ satisfies the following universal property: 
{\it For any Bernstein algebra $(B, \omega')$  and for any element $b\in B$ with $\omega'(b)=1$, there exists a unique baric  homomorphism $\varphi: A \rightarrow B$ such that $\varphi(a)=b$.}\\
Indeed, let $\varphi: A \rightarrow B$ be the linear mapping  defined by $\varphi(a^k)=b^k$ for all $k\geq 1$, equivalently, $\varphi(P(a))=P(b)$ for each polynomial $P\in K[X]$ with no constant term. Then $\varphi$ is baric, that is, $\omega'\circ \varphi=\omega$, since $\omega'(\varphi(P(a))=\omega'(P(b))=P(1)=\omega(P(a))$. In addition, $\varphi$ is multiplicative, since $$\varphi(a.a^k)=\varphi(a^{k+1})=b^{k+1}=b.b^k=
\varphi(a).\varphi(a^k),$$
and, for $k, l\geq 2, \;\; \varphi(a^l.a^k)=\varphi\left(\frac 12(a^l+a^{k})\right)=\frac 12(b^l+b^{k})=b^l.b^k=
\varphi(a^l).\varphi(a^k)$.\\
Furthermore, we have $\text{Im}(\varphi)=\text{alg}(b)$ and $\ker(\varphi)=\{P(a)\; |\; P(b)=0\}$. In particular, as in the situation of Theorem~\ref{th:principal_subalgebra}(ii), if $b$ is algebraic with minimal polynomial $p_b$, then the ideal $\ker(\varphi)=\{P(a)\in A\; |\; p_b |P\}=\{(Qp_b)(a)\; |\; Q\in K[X]\}$ is the subspace spanned by the set $\{L_a^j(p_b(a))\; |\; j\in \N\}$ (see Remark \ref{minimal polynomial}), and  $A/\ker(\varphi) \cong \text{alg}(b)$. It is noteworthy that the ideal $\ker(\varphi)$ is baric, that is $\ker(\varphi) \subseteq \ker(\omega)$, which comes from the fact that $\omega(p_b(a))=p_b(1)=\omega'(p_b(b))=\omega'(0)=0$.\\
As a consequence, the Bernstein algebra in Theorem~\ref{th:principal_subalgebra}(ii) is isomorphic to a quotient algebra of the one in Theorem~\ref{th:principal_subalgebra}(iii).
}
\end{Note}

Next, we shall apply the previous theorem in order to detect some explicit forms of the minimal polynomial $p_a$ of any algebraic element in a given Bernstein algebra $A$. 
	
	\begin{Proposition}
	\label{prop:min_pols_for_algebraic_elements}
	    Let $A$ be a Bernstein algebra, and let $a\in A$ be algebraic with minimal polynomial $p_a$. Then:\\
	  (i)  $\deg p_a=1 \iff p_a(X)=X \iff a=0$;\\
	  (ii) $\deg p_a=2 \iff p_a(X)=X^2-\omega(a)X$;\\
	  (iii) $\deg p_a=3 \iff p_a(X)=X^3-\omega(a)X^2$;\\
	  (iv) $\deg p_a \geq 4 \iff p_a(X)=(X^3-\omega(a)X^2)Q(X)$ for some polynomial $Q(X)$ with $\deg Q\geq 1$.
	    	\end{Proposition}

	\pr (i) is trivial.\\
	(ii) Assume that $a\neq 0$ is algebraic of degree 1, that is, $\deg p_a=2$. Thus, $a^2=\lambda a$ for some $\lambda \in K$. This implies $\omega(a)^2=\lambda \omega(a)$, and so $\lambda =\omega(a)$ whenever $\omega(a)\neq 0$. Now, if $\omega(a)=0$, then the identity $(a^2)^2=\omega(a)^2a^2$ yields $\lambda ^3a=0$, and so $\lambda =0$. It follows that $a^2=\omega(a)a$ in either case, which shows (ii).   \\
	(iii) The condition $\deg p_a=3$ means $\dim \text{alg}(a)=2$. According to Proposition \ref{Conseibo}(i), $\text{alg}(a)$ is isomorphic to the constant Bernstein algebra $A_1$ of type $(1, 1)$. This Bernstein algebra, also being  Jordan, satisfies the equation $x^3-\omega(x)x^2=0$.  By  (ii), we have $p_a(X)\neq X^2-\omega(a)X=0$. Therefore, $p_a(X)=X^3-\omega(a)X^2$. \\
	(iv) Suppose now that $\deg p_a\geq 4$, i.e., $\dim \text{alg}(a)\geq 3$. Let us start by the case $\dim \text{alg}(a)\geq 4$ and write $\text{alg}(a)=Ke_a\oplus U_a\oplus V_a$ as in Theorem \ref{th:principal_subalgebra}(ii), where $\dim U_a\geq 2$. Applying identity \eqref{eq:LHS_train_pol_via_Lv} to $a=e_a+2u_1+v_1$,  we have 
	$$
	\left(a^3-\omega(a)a^2  \right) \left(a-\frac 12\omega(a)\right)^k=L_{v_1}^k(a^3-\omega(a)a^2), \quad \mbox{ for all }k\geq 1.
	$$
	On the other hand, since $\text{alg}(a)$ is finite-dimensional, the restriction $l_{v_1}$ of the operator $L_{v_1}$ on $U_a$ is algebraic. Therefore, as $a^3-\omega(a)a^2\in U_a$, we have $\sum\limits_{k=0}^{n-2} \alpha_k L_{v_1}^{k}(a^3-\omega(a)a^2)=0$ for some scalars $\alpha_0, \dots, \alpha_{n-2}$, not all zero. It follows that 
	$$\left(a^3-\omega(a)a^2  \right)\sum\limits_{k=0}^{n-2}\alpha_k\left(a-\frac 12\omega(a)\right)^k=0.$$
	Consequently,  $a$ is annihilated by the non zero polynomial $T(X)=\left(X^3-\omega(a)X^2  \right)Q(X)$, where $Q(X)=\sum\limits_{k=0}^{n-2}\alpha_k\left(X-\frac 12\omega(a)\right)^k$. \\
	It remains the case $\dim \text{alg}(a)=3$ where we utilise Proposition \ref{Conseibo}(ii). Applying identity (2.5) to $a=e+2u_1+v_1$ for $k=1$, we get 
	$$
	\left(a^3-\omega(a)a^2  \right) \left(a-\frac 12\omega(a)\right)=L_{v_1}(a^3-\omega(a)a^2)=v_1u_1=(\alpha-\frac 32)u_1=(\alpha-\frac 32)(a^3-\omega(a)a^2),$$
	deducing that $$(a^3-\omega(a)a^2)[a-(\alpha-1)]=0$$
	Thus, the polynomial $T(X)=(X^3-\omega(a)X^2)[X-(\alpha-1)]$ annihilates $a$. 	\ep\\
	
	As announced in Remark \ref{minimal polynomial}, the preceding proposition implies that if $A$ is a Bernstein algebra and if $a\in A$ is algebraic of degree $n\geq 2$ satisfying the equation 
	$$a^{n+1}=\gamma_1 a+ \gamma_2 a^2 +\dots+ \gamma_{n}a^n,$$
		then $\gamma_1=0$. \\

\section{Bernstein algebras that are locally train}

\label{sec:locally_train}
\label{sec:train}

As explained in Section \ref{sec:alg_Bernstein_algebras}, a Bernstein algebra $A=Ke\oplus U\oplus V$ is algebraic if and only if the operators $L_v: U \rightarrow U$ are locally algebraic (cf. Theorem~\ref{th:A_is_alg_iff_Lv_is_alg}). An important particular case of locally algebraic operators is the notion of locally nilpotent operators. An endomorphism $f$ of a vector space $X$ is said to be {\it locally nilpotent} \cite[page 38]{Kap} if every $x\in X$ is annihilated by some power of $f$, that is, $f^{n_x}(x)=0$ for some integer $n_x$ depending on $x$. Now, the natural question is to know what
sort of algebraic Bernstein algebras $A=Ke\oplus U\oplus V$ have the property that all the operators $L_v: U \rightarrow U$ are locally nilpotent. In order to answer this question, we remind from \cite[page 1164]{Ouatt1} that a baric algebra $(A, \omega)$ is called a {\it locally train algebra} if the subalgebra
alg$(x)$ is a train algebra for every $x\in H$. Clearly, this is equivalent to saying that the subalgebra
alg$(x)$ is a train algebra for every $x\in A$ with $\omega(x)\neq 0$. 
\begin{Note}
\label{ex:locally_train_non_train}
    {\rm Evidently, train algebras are locally train, but the converse is generally false, at least in the class of Bernstein algebras. Indeed, let $A$ be the Bernstein algebra considered in Example~\ref{ex:algebraic_unbounded_degree}.
    Let $y=\sum\limits_{i=1}^n \beta_i u_i+\lambda v$ be an arbitrary element 
    of $N$ in this example. Then $y^{n+1} =0$, hence $N$ is nil. It follows that in each subalgebra
    alg$(x)$ the barideal $N(\alg(x))$ is nil,  so $A$ is locally train by Theorem~\ref{th:f_dim_train}. On the other hand, although $A$ is algebraic, it is not train because of its unbounded degree.
    }
\end{Note}



We start by stating a property of locally train Bernstein algebras that provides further examples of algebraic Bernstein algebras:

\begin{Proposition} 
\label{prop:locally_train_iff_fdim_subalgebra_is_train}
For a Bernstein algebra $A$, the following conditions are equivalent:\\
(i) $A$ is locally train\\
(ii) $A$ is algebraic and every finite-dimensional baric subalgebra of $A$ is a train algebra.
\end{Proposition}

\pr (i) $\Rightarrow$ (ii): Let $x\in A$ with $\omega(x)\neq 0$. Then the subalgebra alg$(x)$ is train, and so finite-dimensional by 
Theorem~\ref{th:principal_subalgebra}. Therefore $A$ is algebraic in view of Proposition \ref{alg H}.
 
 Let now $B$ be a finite-dimensional baric subalgebra of $A$. Since $B$ is a finite-dimensional Bernstein algebra such that each alg$(x)$ for $\omega(x)\neq 0$ is train, 
 it follows from \cite[Proposition 4.1]{Cons} that $B$ is a train algebra.\\
 (ii) $\Rightarrow$ (i): Conversely, under the assumptions of (ii), any element $x\in A$ with $\omega(x)\neq 0$ generates a finite-dimensional subalgebra $\text{alg}(x)$, which is then a train algebra. Hence, $A$ is locally train.
 \ep\\

The previous proposition yields directly:
\begin{Corollary} 
\label{cor:fdim_loc_train_is_train}
Let $A$ be a finite-dimensional Bernstein algebra. If $A$ is locally train, then it is train.
\end{Corollary}

The following technical formula, mentioned by Bayara in \cite[page 364]{Bayara}, will be useful for our applications.  The proof is straightforward and omitted.  
\begin{Lemma}
\label{lem:L^{k+3}_{x|N}}
 In any Bernstein algebra $A=Ke\oplus U\oplus V$,  the following identity holds:
 $$L^{k+3}_{x|N} =L^k_{v|U}\circ L^3_{x|N} \mbox{  for all } x=u+v\in N \mbox{ and } k\geq 0.$$
\end{Lemma}

In \cite{Wa1} some importance  was given to the following polynomials in arbitrary Bernstein algebras:
$$f_k(x)=(x^3-\omega(x)x^2)\left(x-\frac 12 \omega(x)\right)^{k-3}\; \mbox{for }k\geq 3,$$
satisfying the recursive relation 
$$f_{k+1}(x)=xf_k(x)-\frac 12 \omega(x)f_k(x) \;\mbox{ for all }k\geq 3.$$
In particular, we have pointed out in Theorem \ref{th:f_dim_train} that if a finite-dimensional Bernstein algebra $A$ is a train algebra of rank $r\geq 3$, then it satisfies the train equation
$f_r(x)=0$. We will later extend this result to the infinite-dimensional context. 

\begin{Lemma}
\label{lem:Lv_is_k-nilp}
Suppose that, for each $v\in V$, the restriction $l_v$ of the operator $L_v$   
to the Lyubich ideal $L(A)$ is nilpotent. Then for each $x\in A$ the equality $f_m(x)=0$
holds for some $m$. Moreover, if $l_v$ is $k$-nilpotent for each $v\in V$, then 
the identity $f_{k+3}(x)=0$ holds in $A$.
\end{Lemma}

\pr
Similarly to the proof of Proposition~\ref{pr:2.18}, let $x=\alpha e+u+v\in A$. Then for each $k\ge 1$ we have
$$
f_{k+3}(x) = \left(x-\frac 12 \omega(x)\right)^{k} (x^3-\omega(x)x^2) = L_v^k(x^3-\omega(x)x^2),
$$
where $x^3-\omega(x)x^2 \in L(A)$. If $l_v^k =0$, then $f_{k+3}(x) = l_v^k(x^3-\omega(x)x^2) = 0$.
\ep\\

Let $A$ be an arbitrary Bernstein algebra and let $a\in A$ with $\omega(a)\neq 0$. We will say that $a$ is a {\it train element} if it satisfies $f_k(a)=0$ for some $k\geq 3$. The next result states that a sufficient condition for the Bernstein subalgebra $\text{alg}(a)$ to be train is that the generator $a$ be a train element. 

\begin{Proposition}
\label{pr:train_element}
Let $A$ be an arbitrary Bernstein algebra and let $a\in A$ with $\omega(a)\neq 0$. If $a$ is a train element satisfying $f_m(a)=0$ for some minimal integer $m\geq 3$, then: \\
(i) $\;a$ is algebraic with minimal polynomial 
$f_m(X)$;\\
(ii)  $\;\alg (a)$ is a train algebra of rank $m$ and train equation $f_m(x)=0$.
\end{Proposition}
\pr
(i) As $\omega(a)\neq 0$, there is no loss of generality to assume $\omega(a)=1$. In this case $a$ satisfies the polynomial equation $f_m(a)=(a^3-a^2)\left(a-\frac 12\right)^{m-3}=0$, so  $a$ is annihilated by the polynomial $f_m(X)=(X^3-X^2)\left(X-\frac 12\right)^{m-3}$ of degree $m$, which shows that $a$ is algebraic. Let $p_a(X)\in K[X]$  be the minimal polynomial of $a$ defined as in Remark \ref{minimal polynomial}. Then $p_a(X)$ divides the polynomial $f_m(X)$. But $p_a(X)$ is divisible by $(X^3-X^2)$ in view of  Proposition~\ref{prop:min_pols_for_algebraic_elements}.
It follows that $p_a(X)=(X^3-X^2)\left(X-\frac 12\right)^{k} = f_{k+3}(X)$ for some $k$ with $0\leq k\leq m-3$.
By minimality of $m$, we get $k+3=m$.

(ii) Now, we  show that $\text{alg}(a)$ is in fact a train algebra. Consider the Peirce decomposition $\text{alg}(a)=Ke_a\oplus U_a\oplus V_a$ of the Bernstein algebra $\text{alg}(a)$ with respect to the idempotent $e_a=a^2$. According to Theorem~\ref{th:principal_subalgebra}, we have $U_a^2=0$ and $V_a=<v_1>$ with $v_1=a+a^2-2a^3$. Furthermore, since $a$ is algebraic of rank $n=m-1$, there is a basis $(u_1, \dots, u_{n-2})$ of $U_a$ with $u_1=a^3-a^2$, satisfying the following multiplication table:
$$ v_1u_i=u_{i+1}\; (1\leq i\leq n-3), \;\;\;\; v_1u_{n-2}=\beta_1u_1+\dots+\beta_{n-2}u_{n-2},\;\;\;\; v_1^2=-2u_1-4u_2$$
In virtue of Theorem~\ref{th:principal_subalgebra}(ii), the subalgebra $\text{alg}(a)$ will be train if one establishes that all $\beta_i$ are zero. To this end, we shall exploit the equation $f_m(a)=0$ by observing that $a=e_a+2u_1+v_1$. For this, applying identity \eqref{eq:LHS_train_pol_via_Lv} to $a=e_a+2u_1+v_1$ for $k=n-2$,  we have 
$$	(a^3-a^2)\left(a-\frac 12\right)^{n-2}=L_{v_1}^{n-2}(a^3-\omega(a)a^2)=L_{v_1}^{n-2}(u_1),$$
that is,
$$f_{n+1}(a)=\beta_1u_1+\dots+\beta_{n-2}u_{n-2}.
	$$
Taking into account the equation $f_m(a)=0$, where $m=n+1$, it follows from the latter relation that all $\beta_i$ vanish. We conclude from Theorem \ref{th:f_dim_train} that $\text{alg}(a)$ is a train algebra.  \ep\\

%
%
%

The next result characterizes locally train Bernstein algebras in connection with locally nilpotent operators and offers a natural generalization of Theorem~\ref{th:f_dim_train}.
 \begin{Theorem} 
 \label{th:locally_train}
 Let $A=Ke\oplus U\oplus V$ be a Bernstein algebra. The following
conditions are equivalent:\\
(i) $A$ is a locally train algebra;\\
(ii) For all $v\in V$, the operators $L_v: U \rightarrow U$ are locally nilpotent. \\
(iii) For all $v\in V$, the operators $L_v: L(A)  \rightarrow L(A) $ are locally nilpotent. \\
(iv) For all $x\in N$, the operators $L_x: N \rightarrow N$ are locally
nilpotent.
\end{Theorem}

\pr 
(i) $\Rightarrow$ (ii): Assume $A=Ke\oplus U\oplus V$ is locally train and let $u\in U$ and $v\in V$. Consider the Bernstein subalgebra $B=\text{alg}(e, u, v)$ of $A$ generated by the three elements $e, u, v$. By direct  computations involving
\cite[Lemma 3.1]{GarciaGonzalez} one may  show that
$$B=\langle e,\; L_v^i(u)\ (i\geq 0), \; L_v^j(v)\ (j\geq 1), \; L_v^k(u^2)\ (k\geq 0)\rangle  $$
$$\qquad=\;\langle e,\; L_v^i(u)\ (i\geq 0), \; L_v^i(v^2)\ (i\geq 0), \;L_v^i(u^2)\ (i\geq 0)\rangle ,$$ 
where other products are zero. 
Now, since $A$ is locally train, then $A$ is algebraic according to Proposition~\ref{prop:locally_train_iff_fdim_subalgebra_is_train}.
We infer from Theorem~\ref{th:A_is_alg_iff_Lv_is_alg} that the operator $L_v: U \rightarrow U$ is locally algebraic. As a consequence, the subspaces $\langle L_v^i(u) \
(i\geq 0)\rangle , \; \langle L_v^i(v^2) \ (i\geq 0)\rangle  $ and $\langle L_v^i(u^2) \
(i\geq 0)\rangle = \langle u^2, \; L_v^i(vu^2) \ (i\geq 1)\rangle $ of $A$ are
finite-dimensional, and so also is $B$. On the other hand, by hypothesis and 
Corollary~\ref{cor:fdim_loc_train_is_train}
(see also \cite[Proposition 4.1]{Cons}), the finite-dimensional Bernstein algebra $B$ is a train algebra of some rank. 
It follows from  Theorem~\ref{th:f_dim_train}
that 
the operator  $L_v: U_B \rightarrow U_B$ is nilpotent. Finally, since $u\in U_B$, we conclude that $L_v^n(u)=0$ for some positive integer $n$.

(ii) $\Rightarrow$ (iv): It is enough to apply 
Lemma~\ref{lem:L^{k+3}_{x|N}}.

(iii) $\Rightarrow$ (i): By Lemma~\ref{lem:Lv_is_k-nilp}, each element $x$ of $A$ satisfies a train equation of some rank. Then the algebra $C = \alg(x)$ is algebraic, hence finite-dimensional. 
 Applying again Lemma~\ref{lem:Lv_is_k-nilp} to each element $y\in N(C) = N(A)\cap C$, we get $y^{k_y+3}=0$ for some integer $k_y$.
Then all elements of $N(C)$ are nilpotent, and  Theorem~\ref{th:f_dim_train} implies that $C$ is train. Thus, $A$ is locally train.

Finally, the proof of the theorem is achieved, since the implication (iv) $\Rightarrow$ (iii) is obvious.\ep \\


Now, we study the behavior of the barideal $N=\ker(\omega)$ of a locally train Bernstein  algebra. We know from Theorem~\ref{th:A_is_alg_iff_Lv_is_alg} 
that $N$ is algebraic, but we want to give a
more precise information on $N$ in this case. The next result provides the answer to this question.

\begin{Theorem} 
\label{th:N_is_nil}
Let $A=Ke\oplus U\oplus V$ be a Bernstein
algebra and let $N=\ker(\omega)=U\oplus V$. Then $A$ is locally train if and only if the ideal $N$ is nil.
\end{Theorem}

\pr 
If $A$ is locally train, then by Theorem~\ref{th:locally_train} (iv) the operators 
$L_x:N\to N$ are nilpotent for each $x\in N$, so that $x^n = L_x^{n-1}(x) = 0$ for $n$ large enough.
Thus, $N$ is a nil-algebra.

Reversely, if $N$ is a nil-algebra, then $N$ is algebraic. Indeed, every $x\in N$ satisfies $x^{n_x}=0$ for some $n_x$,
so alg$(x)=\langle x,x^2, \dots, x^{n_x-1}\rangle $ is finite-dimensional (since $x^ix^j=0$ if $i,j\geq 2$). It follows from Theorem~\ref{th:A_is_alg_iff_Lv_is_alg} 
that $A$ is algebraic. Thus, for any $a\in H$, $A_a=\mbox{alg}(a)$ is a finite-dimensional baric subalgebra of $A$. Now, since $\ker(\omega_{|A_a})=N\cap A_a\subset N$, then $\ker(\omega_{A_a})$ is a fortiori a nil-algebra. Applying Theorem~\ref{th:f_dim_train} to  $A_a$, 
we infer that $A_a$ is a train algebra. This proves that $A$ is locally train. \ep \\

The next theorem is an extension of Theorem~\ref{th:f_dim_train} to the infinite-dimensional case. Before doing this, 
we use the following result obtained by Bayara in \cite[Th\'eor\`eme 3.3.4]{Bayara-thesis}. Recall  that
the identity $(x^2)^2=0$ holds for each $x\in N(A)$ in a Bernstein algebra $A$. The proof of Bayara is itself adapted from the proof of Walcher \cite[Page 161]{Wa1} for Bernstein train algebras. 

\begin{Lemma}
\label{lem:Bayara_on_nil}
Suppose that the identity $(x^2)^2=0$ holds in a commutative (non-associative) algebra $S$.
Then the following conditions are equivalent:

(i) $S$ a nil-algebra of bounded nil-index;

(ii)  $S$ is Engel, that is, the operators $L_x: S \rightarrow S$ for $x\in S$ are
nilpotent with bounded nil-indexes.
\end{Lemma}

\begin{Theorem} 
\label{th:train_algebra}
\label{th:3.8}
Let $A$ be a Bernstein algebra. The following
conditions are equivalent:\\
(i) $A$ is a train algebra;\\
(ii) For all $v\in V$, the operators $L_v: U \rightarrow U$ are nilpotent with
bounded nil-indexes. \\
(iii) For all $v\in V$, the operators $L_v: L(A)  \rightarrow L(A) $ are nilpotent with
bounded nil-indexes.\\
(iv) For all $x\in N$, the operators $L_x: N \rightarrow N$ are
nilpotent with bounded nil-indexes.\\
(v) $N$ is a nil-algebra of bounded nil-index.\\
In this case, the train equation of $A$ has the form
$(x^3-\omega(x)x^2)(x-\frac 12\omega(x))^{r-3}=0$, where $r$ is the rank of $A$.
\end{Theorem}

\pr
The implications (iv)$\Rightarrow$(ii)$\Rightarrow$(iii)
are obvious. Since the  train equation~(\ref{train}) of degree $r$ implies the equality $x^r=0$ for 
$x\in N$, the implication (i)$\Rightarrow$(v) is trivial too. 
Moreover, the implication
(iii)$\Rightarrow$(i) follows from 
Lemma~\ref{lem:Lv_is_k-nilp}. 
 The rest implication (v)$\Rightarrow$(iv) is a particular case of Lemma~\ref{lem:Bayara_on_nil}.
 \ep





	\section{Algebraic Banach Bernstein algebras}
	
	\label{sec:banach}
	
	Our objective in this section is to  establish some
	results  about Banach Bernstein algebras that are algebraic. Recall
	\cite{Ly2} that a {\it Banach Bernstein algebra} $(A, \omega, \|.
	\|)$ is a Bernstein algebra $(A, \omega)$ over the field $K=\R$ or
	$\Complex$ equipped with a structure of Banach algebra $(A, \|. \|)$. It was shown  in
	\cite{B} that the weight function $\omega: A \rightarrow K$ is automatically continuous. Hence, the ideal $N=\ker(\omega)$ is closed and is also a Banach algebra.
	Bernstein algebras were studied intensely from the algebraic point
	of view, but Banach Bernstein algebras are much less known.  
	
	 Returning to our investigation, we have seen in Example~\ref{ex:algebraic_unbounded_degree} that the
	hypothesis that $A$ be of bounded degree in Theorem~\ref{th:Theorem_2.15} cannot be
	relaxed. In other words, if $A$ is a Bernstein algebra that is
	algebraic, then the operators $L_v: U \rightarrow U$ are not
	necessarily algebraic. However, when the Bernstein algebra $A$ is
	in addition a Banach algebra, this situation cannot hold. Namely,
	we may delete the ``bounded degree'' hypothesis in Theorem~\ref{th:Theorem_2.15}. 
	
	
	
	
	\begin{Theorem} 
	\label{th:banach1}
	Let $A=Ke\oplus U\oplus V$ be a Banach Bernstein
		algebra. The following assertions are equivalent:\\
		(i) the algebra $A$ is algebraic;\\
		(ii) the operators $L_v: U \rightarrow U$ are algebraic for
		all $v\in V$;\\
		(iii) the algebra $A$ is algebraic of bounded degree;\\
		(iv) the  operators $L_v: U \rightarrow U$ are algebraic and the degrees are bounded collectively for 
		all $v\in V$.
	\end{Theorem}
	
	{\it Proof. \ } In view of  Theorems~\ref{th:Lv_is_alg_bounded} and~\ref{th:A_is_alg_iff_Lv_is_alg}, the implication $(ii) \Rightarrow (i)$ and the equivalence $(iv) \Leftrightarrow (iii)$ are already
	true in full generality without the additional Banach hypothesis.\\
	$(i) \Rightarrow (iii)$: Let $v\in V$. In virtue of Theorem~\ref{th:A_is_alg_iff_Lv_is_alg}, the operator $L_v$ is
	locally algebraic. 
	It is known \cite[Corollary 2.8.33]{Rod} that any arbitrary nonassociative Banach algebra which is algebraic is of bounded degree. 
	\ep\\
	
	
	
	
	
	\begin{Theorem}
	\label{th:banach2}
	    Let $A$ be a Banach Bernstein algebra. If $A$ is locally train, then $A$ is train.
	\end{Theorem}
	
	\pr If $A$ is locally train, by Theorem~\ref{th:N_is_nil} the ideal $N=\ker(\omega)$ is nil. It follows from \cite[Corollary 2.8.41]{Rod} that the Banach nil-algebra $N$ is of bounded nil-index. The proof is concluded by applying Theorem~\ref{th:train_algebra}. \ep\\
	
	This result permits us to recover the finite-dimensional case in Theorem~\ref{th:f_dim_train}. 
	
	\section{Algebraic Bernstein algebras with low degrees}
	
		\label{sec:alg_low_degree}

	In this section we deal with algebraic Bernstein algebras with low bounded degrees.
	First, we begin with the case of a  Bernstein algebra which is
	algebraic with bounded degree 1, that is $x$ and $x^2$ are
	linearly dependent for all $x\in A$. The following fact is actually valid for arbitrary commutative baric algebras.
	
	\begin{Proposition} 
	\label{prop:former4.1}
	\label{prop:baric_bound_deg1}
	Let $A$ be an arbitrary commutative baric algebra which is algebraic of bounded degree $1$. Then $A$ is a Bernstein algebra satisfying the train equation $x^2=\omega(x)x$.
	\end{Proposition}
	
	\pr By hypothesis, for each $x\in A$ there
	exists $\lambda_x\in K$ such that $x^2=\lambda_x x$, which gives
	$\omega(x)^2=\lambda_x\omega(x)$. This implies
	$\lambda_x=\omega(x)$ and so $x^2=\omega(x)x$ whenever
	$\omega(x)\neq 0$. Now, since the set $\{x\in A \; | \;
	\omega(x)\neq 0\}$ is dense for the Zariski topology, it follows
	that $x^2=\omega(x)x$ for all $x\in A$. Hence, $A$ is an
	elementary Bernstein algebra. \ep\\
	
	
	\quad Now, let $A$ be an arbitrary commutative algebra (not necessarily Bernstein) which is algebraic of bounded degree 2. This means that $x, x^2, x^3$ are linearly dependent for all $x\in A$. Equivalently,  every element of $A$ generates a subalgebra of dimension not greater
than two. Such algebras appear in the papers \cite{Wa2, Wa3} by Walcher as {\it algebras of rank three}. Walcher developed the
	structure of these algebras and showed the
	existence of a linear form $\gamma_1: A \rightarrow K$ and a
	quadratic form $\gamma_2: A \rightarrow K$ such that
	$x^3=\gamma_1(x)x^2+\gamma_2(x)x$
	for all $x\in A$.\\
	Now, assume that $A$ is a Bernstein algebra of degree 2, that is $A$ is an algebra of rank 3 in the sense of Walcher. There exists  an element $a$ in $A$ such that $a$ and $a^2$ are linearly independent; otherwise  $\deg(A)=1$. Applying Proposition \ref{prop:min_pols_for_algebraic_elements}(iii), we have $p_a(X)=X^3-\omega(a)X^2$ and $p_y(X)=X^2-\omega(y)X$ whenever $y\in A$ with $\deg(y)=1$. Therefore, $A$ satisfies the identity $x^3-\omega(x)x^2=0$, so $A$ is a Bernstein- Jordan algebra (which does not satisfy the identity $x^2=\omega(x)x$).\\
	We have thus obtained the following result.  
	
	\begin{Theorem} 
		\label{th:former_4.2}
		\label{th:jordan}
		Let $(A,\omega)$ be a Bernstein algebra. The following conditions are 	equivalent:\\
		{\rm (i) } $\; A$ is algebraic of degree $2$.\\
		{\rm (ii) } $\; A$ is a Jordan algebra (that is not an elementary Bernstein algebra).
	\end{Theorem}
	

	
	Applying Proposition~\ref{prop:former4.1} together with Theorem~\ref{th:former_4.2}, 
	we obtain the following consequence.
	\begin{Corollary} 
	\label{cor:former4.4}
	Let $(A,\omega)$ be a Bernstein algebra. If $A$ is algebraic of degree $\leq 2 $, then  $A$ is  a Jordan algebra.
	\end{Corollary}

	\quad Now we give an example of an algebraic Bernstein algebra of bounded degree 3 which is
	not Jordan. This means that the bound 2 in Corollary~\ref{cor:former4.4} is the best possible.
	
	\begin{Example} {\rm This construction is based on the Zhevlakov example \cite[Page 82]{russe} of a (special) Jordan algebra satisfying the identity $x^3=0$, which is solvable but not nilpotent. Let $X=\{x_1, x_2, \dots,\}$ be a countable set of symbols. A word $x_{i_1}.x_{i_2}.\dots x_{i_n}$ in the alphabet $X$ is called regular if $i_1<i_2<\dots < i_n$. Regular words form the basis of the commutative algebra $N$ whose multiplication is given on the basis of:\\
			\quad $\;\;\bullet\;$ $u*v=v*u\;\;$ for all regular words $u, v$ ;\\
			$\;\;\bullet\;$ $x_i*x_j=x_ix_j \;\;\text{if } i<j$ ;\\
			$\;\;\bullet\;$ $(x_{j_1}\dots x_{j_n})*x_i=(-1)^\sigma x_{j_1}\dots x_j \dots x_{j_n}\;\;$ if  $n\geq 2$ and $i>j_1$, where $\sigma$ is the number of inversions in the permutation $(j_1, \dots, j_n, i)$ ;\\
			$\;\;\bullet\;$ other products are zero. \\
			Let us consider the subspace $V$ (respectively, $U$) spanned by the words of degrees 1  (respectively, of degrees $\geq 2$). Then $N=U\oplus V$ and $U^2=0, \, UV\subset U, \,V^2\subset U$. To make the Zhevlakov algebra $N$ into a Bernstein algebra,  we embed $N$ in a
			commutative algebra $A=Ke\oplus N$ by adjoining an
			idempotent $e$ and completing the multiplication table by
			$eu=\frac 12 u$ and $ev=0$ for all $u\in U$ and $v\in V$. Clearly,
			the linear mapping $\omega: A \rightarrow K$ defined by
			$\omega(\alpha e+u+v)=\alpha$ is a weight function. Further, by Lemma~\ref{lemma:Lemma1.1} $(A,
			\omega)$ is actually a Bernstein algebra with Peirce components
			$U$ and $V$ with respect to the idempotent $e$. \\
			On the other hand, since $x^3=0$ for all $x\in N$, it follows from Theorem~\ref{th:train_algebra} that $A$ is a train algebra. In fact, by straightforward calculations we may prove that $A$ is a train algebra of rank 4 satisfying the train equation $(x^3-\omega(x)x^2)(x-\frac 12\omega(x))=0$, or equivalently, $x^4-\frac 32\omega(x)x^3+\frac 12\omega(x)^2x^2=0$.
	}\end{Example}
	
	\begin{Note} {\rm In addition to the infinite-dimensional previous counter-example, it is possible to give an alternative counter-example that is finite-dimensional. 
	It is sufficient to reconsider  the Bernstein algebra $A$ treated in Example~\ref{ex:2.8}, 
	that is algebraic of degree 3. Since $A$ is not train, it is not Jordan at all. This example shows that in general, an algebraic Bernstein algebra of degree $\geq 3$ need not be train.}
	\end{Note}
	
	\section{The  Kurosh problem}
		\label{sec:Kurosh}

	\esp Our purpose in this section is to formulate the Kurosh problem
	for Bernstein algebras in terms of associative algebras. In
	this way, the problem will be equivalent to a natural
	question for associative algebras, which may be
	viewed as a weak form of the classical  Kurosh problem
	for associative algebras.
	
The main goal of this section is to prove the following result.
    
    	\begin{Theorem} 
    	\label{th:Kurosh}
		(a) Each finitely generated Bernstein algebra which is
		algebraic of bounded degree $\leq 2$ is  finite-dimensional.
		
		(b) For each $d\ge 3$ and $n\ge 2$, there exists an  infinite-dimensional  $n$-generated Bernstein algebra which is
		algebraic of bounded degree $d$.
		
		(c) Under the conditions of (b), there exists  an infinite-dimensional finitely $n$-generated Bernstein train algebra of rank $d+1$.
	\end{Theorem}

	\pr 
    (a) Let $A$ be a Bernstein algebra as in Part (a).     	
 By Corollary~\ref{cor:former4.4}, $A$ is also a Jordan algebra. 
	Hence, finite-dimensionality follows from \cite{Su} (see Corollary~\ref{cor:fdim_Jordan} below).
	An explicit minimal example of an  infinite-dimensional  $2$-generated train algebra of rank 4 (which is therefore algebraic of degree 3) is discussed in Lemma~\ref{Lemma:nil_assoc_algebras}.
		The proof of Part (c) consists of Example~\ref{Example: Kurosh for Bernstein} and Remark~\ref{rem:Kurosh_example_n_k}
	given below, while  Part (b) is an immediate corollary of Part (c).\\

	The plan of this section is the following. First, we investigate several cases in which the Kurosh problem for Bernstein algebras have the positive solution. Next, we give a construction of a Bernstein algebra by means of an arbitrary associative algebra in 
	Theorem~\ref{th:Bernstein_via_assoc}. We utilise this construction to prove the results needed in the rest of  Theorem~\ref{th:Kurosh}. \\

	

	Before starting, let us review the deep link
exhibited in \cite{B} between Bernstein algebras and modules over associative (noncommutative) algebras.
Let $A=Ke\oplus U\oplus V$
	be a Bernstein algebra, and consider the free unital associative
	(noncommutative) algebra $T(V)$ 
	generated by the 
	vector space
	$V$. Then  the Lyubich ideal $L(A) $ becomes a left module over
	$T(V)$ by setting
	\begin{equation}\label{product-K<V>}
	    (v_1*\ldots *v_k).u=v_1(\ldots (v_ku)\ldots ), \mbox{ for all }
	v_1, \ldots, v_k\in V \mbox{ and } u\in L(A).
	\end{equation}
		This $T(V)$-module $L(A) $ contains some finiteness information about the 
	Bernstein algebra $A$. Especially, by \cite[Theorem 3.5]{B}  $A$ is finitely generated
	if and only if $A/L(A) $ is finite-dimensional and the
	$T(V)$-module $L(A) $ is finitely generated. In this case, $V$ is a fortiori 
	finite-dimensional, say $V=\langle v_1, \dots, v_n \rangle $, and $L(A) =\sum\limits_{i=1}^r T(V).u_i$ for some elements $u_1, \dots, u_r\in L(A) $.
Therefore, it is a simple matter to check from (\ref{product-K<V>}) that
	\begin{equation} \label{alg}
	L(A) =\sum\limits_{i=1}^r \text{alg}(L_{v_1}, \ldots, L_{v_n}).u_i,
	\end{equation}
	where $\text{alg}(L_{v_1}, \ldots, L_{v_n})$ is the subalgebra of $\text{End}_K(U)$ generated by $L_{v_1}, \ldots, L_{v_n}$ and $\text{alg}(L_{v_1}, \ldots, L_{v_n}).u_i$ means the subspace of $A$ spanned by all elements $T(u_i)$ with $T\in \text{alg}(L_{v_1}, \ldots, L_{v_n})$.\\
	
		\quad Now, we shall prove that the Kurosh problem for
	Bernstein algebras have a positive answer for various special	cases. We start with the following result concerning the case $\dim V=1$.  It turns out that  a  finitely generated Bernstein algebra with $\dim V =1$ that is algebraic will actually be of bounded degree. More specifically: 
	
	\begin{Proposition} 
	\label{prop:dimV=1} Let $A=Ke\oplus U\oplus V$ be a finitely generated Bernstein
		algebra with  $\dim V=1$. If $A$ is algebraic,
		then  $A$ is finite-dimensional.
	\end{Proposition}
	
	{\it Proof. } 
Write $V=\langle v\rangle $, where $v\in V$. Since $A$ is algebraic, then the operator $L_v$ is locally algebraic thanks to Theorem~\ref{th:A_is_alg_iff_Lv_is_alg}.
 Hence, for any $u\in U$ there exists
a nonzero polynomial $P_u(X)=\sum\limits_{i=0}^p \alpha_i X^i$ in $K[X]$ such that $P_u(L_v)(u)=0$. Therefore, $\;
\sum\limits_{i=0}^p \alpha_i L_v^i(u)=0$, or, equivalently, $g.u=0$, where $g=\alpha_0 1_{K\langle
V\rangle}+\alpha_1 v+ \alpha_2 (v*v)+\ldots + \alpha_p(v*\ldots
*v)$ is a nonzero element in $T(V) = K[v]$.
Thus, $L(A)$
 is a torsion module over $T(V)$. It follows from  \cite[Proposition 3.9]{B} that $A$ is finite-dimensional.	 \ep\\
	
	Combining Propositions \ref{prop:locally_train_iff_fdim_subalgebra_is_train} and \ref{prop:dimV=1}, together with Corollary \ref{cor:fdim_loc_train_is_train}, allows us to state:
	
	\begin{Corollary}
		\label{cor:dimV=1} Let $A=Ke\oplus U\oplus V$ be a finitely generated Bernstein
		algebra with \\ $\dim V=1$. If $A$ is locally train, then  $A$ is a finite-dimensional train algebra.
	\end{Corollary}
	
The requirement that $A$ be finitely generated cannot be suppressed in Proposition \ref{prop:dimV=1} and Corollary \ref{cor:dimV=1}. To be convinced, it suffices to return to the infinite-dimensional Bernstein algebra presented in Example \ref{ex:algebraic_unbounded_degree} which is algebraic and locally train (see Remark \ref{ex:locally_train_non_train}), but not finitely generated.

On the other hand, in contrast to the situation $\dim V=1$ in the last two results, it will be shown in Example \ref{Example: Kurosh for Bernstein} that the analogous results are no longer valid when $\dim V\geq 2$.\\

	\quad We have seen in 
	Example \ref{ex:Bernstein_train_is_algebraic}
	that Bernstein-train algebras are special instances of algebraic Bernstein algebras with bounded degrees. 
	Hence, the Kurosh problem in this case asks whether a finitely generated
	Bernstein-train algebra is finite-dimensional. \\On the other hand,  Theorem \ref{th:train_algebra} asserts that a Bernstein algebra $A=Ke\oplus U\oplus V$ is train if and only if the operators $L_v: U \rightarrow U$ $(v\in V)$ are nilpotent of bounded nil-indexes.
	In this approach,
	the next result treats the particular situation where  $L_v^2=0$ for
	all $v\in V$.
	
	
	

	\begin{Proposition} 
	\label{prop:former5.6}
	Let $A=Ke\oplus U\oplus V$ be a finitely generated Bernstein algebra satisfying $L_v^2=0$ for all $v\in V$. Then $A$ is
		finite-dimensional and $N=U\oplus V$ is nilpotent.
	\end{Proposition}
	
	{\it Proof. } 	Making use of equality  (\ref{alg}), we have
	$L(A) =\sum\limits_{i=1}^r \text{alg}(L_{v_1}, \ldots, L_{v_n}).u_i$,
for some  $v_1, \dots , v_n\in V $ and $u_1, \dots, u_r\in L(A)$.
	By hypothesis,  $L_vL_{v'}+L_{v'}L_v=0$ for
	all $v, v'\in V$. This entails that
	$$\text{alg}(L_{v_1}, \ldots, L_{v_n})=\langle L_{v_{i_1}} \ldots L_{v_{i_k}} \;
	| \; 1\leq k\leq n, \; 1\leq i_1 < \ldots< i_k\leq n\}.$$
	Thus, $\text{alg}(L_{v_1}, \ldots, L_{v_n})$ is finite-dimensional,
	and so also is $L(A) =\sum\limits_{i=1}^r \text{alg}(L_{v_1}, \ldots,
	L_{v_n}).u_i$. As $A/L(A)$ is finite-dimensional, then so is $A$. Finally,  the nilpotency of $N$ comes from \cite[Th\'eor\`eme 2.3]{Bayara}.\ep

	\begin{Note}{\rm Let $A=Ke\oplus U\oplus V$ be a Bernstein algebra satisfying $L_v^2=0$. Then  identity  (\ref{eq:LHS_train_pol_via_Lv}) for $k=2$ shows that $A$ satisfies the train equation 
	$$\left( x^3-\omega(x)x^2  \right) \left(x-\frac 12\omega(x)\right)^2=0$$
	In particular, taking $\omega(x)=0$ yields $x^5=0$ for all elements $x\in N=U\oplus V$.

		}
	\end{Note}
	
	\quad The condition $L_v^2=0$ in Proposition~\ref{prop:former5.6} is equivalent to
	$(uv)v=0$ for all $u\in U$ and $v\in V$. In this spirit, recall
	that a Bernstein algebra $A=Ke\oplus U\oplus V$ is a Jordan
	algebra if and only if $(uv)v=0$ and $v^2=0$ for all $u\in U$ and
	$v\in V$. Hence, Proposition~\ref{prop:former5.6} recovers as a direct consequence the
	following well-known result (see \cite{Su}):
	
	
	\begin{Corollary} 
	\label{cor:fdim_Jordan}
	If $A=Ke\oplus U\oplus V$ is  finitely generated Bernstein-Jordan algebra, then $A$ is finite-dimensional and $N=U\oplus V$ is nilpotent.
	\end{Corollary}
	
	Taking into account Proposition \ref{th:jordan}, it follows from the last corollary that any finitely generated Bernstein algebra which is algebraic of
bounded degree $\leq  2$ is  finite-dimensional. Example \ref{Example: Kurosh for Bernstein} below establishes that this result is no longer true for degrees $\geq 3$. \\
	
	\quad In addition to Bernstein-Jordan algebras, the class of nuclear Bernstein algebras also play an important role in the theory of Bernstein algebras. It is worth noting that nuclear Bernstein algebras  are algebraic, since they are train algebras satisfying the train equation $2x^4-3\omega(x)x^3+\omega(x)^2x^2=0$ (Example~\ref{ex:nuclear}). Actually, the Kurosh problem is valid in this class of Bernstein algebras,  since the following result is already well known (see, \cite{Kra, Pe, Su}).

	\begin{Proposition} If $A=Ke\oplus U\oplus V$ is  finitely generated nuclear Bernstein algebra, then $A$ is finite-dimensional and $N=U\oplus V$ is nilpotent.
	\end{Proposition}

	
	\quad In another context, it is known \cite[Corollary 3.6]{B} that each N\oe therian Bernstein algebra is finitely generated, and this result fails for  Artinian algebras. Nonetheless, we can prove the following  Artinian version of Proposition~\ref{prop:former5.6}, which seems to be of independent interest, and whose proof also holds in the finitely generated case:
	
	\begin{Proposition} Let $A=Ke\oplus U\oplus V$ be a Bernstein algebra  satisfying $L_v^2=0$ for all $v\in V$. If $A$ is Artinian, then $A$ is
		finite-dimensional and $N=U\oplus V$ is nilpotent.
	\end{Proposition}
	
	{\it Proof. } Let $I$ be the ideal of $A$ generated by the subspace $V^2$. Then $A/I=K\overline{e}\oplus \overline{U}\oplus\overline{V}$ is a Bernstein Jordan algebra, since it satisfies $\left(\overline{V}\right)^2=0$ and the  identity $(\overline{u}.\overline{v})\overline{v}=0$ (because $(uv)v=L_v^2(u)=0$).
	Since $A$ is Artinian, and so also is $A/I$, it follows from \cite[Theorem 2.3]{B} that $A/I$ is finite-dimensional. \\
	It remains to prove that $I$ is finite-dimensional. Clearly, $I=\sum\limits_{k\geq 2} V^k$, because $eV^k=V^k$ ($V^k\subset U$) and $UV^k=0$ for all  $k\geq 2$.
	According to \cite[Proposition 3.1]{B}, since $A$ is Artinian, then $V$ is finite-dimensional.  Choose a basis $\{v_1, \dots, v_p\}$ of $V$.
	Since $L_v^2(V^2)=0$ and $\left(L_vL_{v'}+L_{v'}L_v\right)(V^2)=0$ for all $v, v'\in V$,
	every monomial $v_{i_1}(v_{i_2}(\ldots v_{i_{k-2}}(v_{i_{k-1}}v_{i_k})\ldots))$
	with $\{i_1, i_2, \dots, i_k\}\subset \{1, \dots, p\}$ vanishes whenever $k>p+2$.
	It follows that $I=V^2+\dots+V^{p+2}$ is finite-dimensional, and so also is $A$.\ep
	\\
	
	In passing, it should be emphasized that the analogs of Proposition \ref{prop:dimV=1} and Corollary \ref{cor:dimV=1}  are false in the Artinian case. Once again, Example~\ref{ex:algebraic_unbounded_degree} is an infinite-dimensional locally train algebra with $\dim V=1$ (see Remark \ref{ex:locally_train_non_train}) which is Artinian by \cite[Example 3.10]{B}.\\
     
	Returning to Proposition \ref{prop:former5.6}, we will see later in Example~\ref{Example: Kurosh for Bernstein} that this result becomes false if  replacing the condition $L_v^2=0\; (v\in V)$ by $L_v^3=0\; (v\in V)$. However, it is valid under the stronger hypothesis $L_x^3=0\; (x\in N)$:

		\begin{Proposition} \label{L_x^3=0}
		Let $A=Ke\oplus U\oplus V$ be a finitely generated Bernstein
		algebra such that the operators $L_x: N \rightarrow N$ satisfy $L_x^3=0$ for all $x\in N=\ker(\omega)$. Then $A$ is
		finite-dimensional and $N$ is nilpotent.
	\end{Proposition}
	
	{\it Proof. } Since $A$ is finitely generated, then also is the subalgebra $N$ \cite[Lemma 7]{Kra}. Thus, $N$ is a finitely generated commutative algebra satisfying the identity $L_x^3=0$. The variety of such algebras $N$  was studied in \cite{Correa} and shown to be nilpotent \cite[Theorem 6]{Correa}. Hence $N$ is nilpotent, and by a straightforward
	argument, $N$ is finite-dimensional. \ep\\
	
	
	 Once again, Example~\ref{Example: Kurosh for Bernstein} below will prove that the above result does not hold when   the condition $L_x^3=0 \;(x\in N)$ is replaced by the weaker condition $L_x^4=0 \; (x\in N)$.\\
	


	
\esp At present, we concentrate our efforts on the Kurosh problem for Bernstein algebras in the general case, and we explain it in the language of associative algebras. Let $A=Ke\oplus U\oplus V$ be a finitely generated Bernstein algebra. Then 
(\ref{alg}) yields 	$L(A) =\sum\limits_{i=1}^r \text{alg}(L_{v_1}, \ldots, L_{v_n}).u_i$
for some $u_1, \dots, u_r\in L(A)$ and $v_1, \dots , v_n\in V $ with $V=\langle v_1, \ldots, v_n\rangle$. If $A$ is algebraic of bounded
degree, then by Theorem \ref{pr:2.18}  the operators $L_v: U \rightarrow U\; (v\in V)$ are all algebraic and
their algebraic degrees are collectively bounded. Consider the finitely generated associative algebra $C=\text{alg}(L_{v_1}, \ldots, L_{v_n})$ and the finite-dimensional vector space $S=\langle L_{v_1}, \ldots, L_{v_n}\rangle$ generating $C$. Then all elements of $S$ are algebraic in $C$ and
their algebraic degrees are uniformly bounded. Besides, as remarked in the proof of Proposition \ref{prop:former5.6},   $A$ is finite-dimensional if and only if $C$ is.\\
Conversely, the next theorem assigns a finitely generated Bernstein algebra to any finitely generated associative algebra $C$.

		\begin{Theorem}
			\label{the:Bernstein_via_assoc}
			\label{th:Bernstein_via_assoc}
			 Suppose that $C$ is a finitely generated associative algebra, and let $S\subset C$ be a finite-dimensional vector space which generates $C$.
		Then  the Cartesian product
		$A = K \times C \times S$ with  coordinate-wise addition and  multiplication defined by
		$$(\alpha, a_1, b_1)(\beta, a_2, b_2)=\left(\alpha\beta,\, \frac 12(\alpha a_2+\beta a_1)+a_1b_2+a_2b_1,\, 0\right)$$
		is a finitely generated exceptional Bernstein algebra.
		
		If, in addition, all elements of $S$ are algebraic of bounded degree in $C$,
		then $A$ is algebraic of bounded degree.
		
		Moreover, if all elements of $S$ are nilpotent of bounded index, then $A$ is a train algebra.
	\end{Theorem}

	{\it Proof. \ } The construction of this algebra $A$ is inspired from \cite[Example]{semiprime} in somewhat altered form (see, also, \cite[page 2577]{B} for an equivalent construction).
Obviously, 
$A$ is a commutative baric algebra with weight function  given by $\omega(\alpha, a, b)=\alpha$. 
The nonzero idempotents of $A$ are all the elements $(1, a, 0)$ with $a\in C$. Choosing the idempotent $e=(1, 0,0)$, we have the  direct sum $A=Ke\oplus U\oplus V$, where $U= (0, C, 0)$ and $V = (0, 0, S)$. Furthermore,  $$eu=\frac 12u \;\, \forall u\in U, \;\, ev=0\; \forall v\in V, \;\, UV\subset U, \;\, U^2=0, \;\, V^2=0$$ It follows from Lemma~\ref{lemma:Lemma1.1} that $A$ is actually an exceptional Bernstein algebra.

	Now, let $S=\langle a_1, \dots, a_n \rangle$ be a span of some elements  $a_1, \dots, a_n$ in $C$.
Since the algebra $C$ is generated by $a_1, \dots, a_n$, then  $A$ is generated by the  set $\{e, (0,a_i,0),  (0,0,a_i)\; | \;1\leq i\leq n\}$.
If $v=(0,0,b)\in V$ with $b\in S$, then $b$ is algebraic by hypothesis, so $P_b(b)=0$ for some nonzero polynomial $P_b\in K[X]$ without a constant term. Let us show that the multiplication operator $L_v: U \rightarrow U$ is algebraic. Indeed, for each $u=(0,a,0)\in U$ with $a\in C$, an easy induction involving associativity in $C$ gives $L^i_v(u)=(0, ab^i,0)=(0, l^i_b(a), 0)$ for all $i\geq 1$, where $l_b: C \rightarrow C$ is the multiplication operator by $b$. Therefore,  $P_b(L_v)(u)=(0, P_b(l_b)(a), 0)=(0,0,0)$, because $P_b(l_b)(a)=P_b(b)a=0$.
Further, it is clear that the degrees of the operators $L_v\; (v\in V)$ are bounded collectively, since the same thing is true for the operators $l_b\; (b\in S)$. Then $A$ is algebraic of bounded degree by Theorem~\ref{th:Lv_is_alg_bounded}.\\
Finally, to prove the last assertion, 
 we apply Theorem~\ref{th:train_algebra} and carry out an analogous reasoning. \ep\\

Actually, in the light of Theorem \ref{th:Bernstein_via_assoc} and the observation preceding it, the Kurosh problem for Bernstein algebras is equivalent to the following question on associative algebras: If $C=\text{alg}(a_1, \dots, a_n)$ is a finitely generated associative algebra such that all elements of the  subspace $S=\langle a_1, \dots, a_n\rangle$ are algebraic of bounded degrees collectively, is $C$ finite-dimensional? The next lemma shows that the latter question has a negative answer in general. 

	\begin{Lemma}
		\label{Lemma:nil_assoc_algebras}
		Let $A_{n,k}$  (where $k,n\ge 1$) be the non-unital associative algebra over an arbitrary field with $n $ generators subject to the following relations: the $k$-th power of each linear combination of the generators is zero. Then $A_{n,k}$ is infinite-dimensional if and only if $n\ge 2$ and $k\ge 3$.
	\end{Lemma}
	
	\pr
	We have $$
	A_{n,k} =  K\!\langle x_1, \dots, x_n |
	(\alpha_1 x_1 +\dots +\alpha_n x_n)^k =0 \mbox{ for all }
	\alpha_1 ,\dots ,\alpha_n \in K \rangle .
	$$
	Obviously, $A_{1,k} =  K\!\langle x_1 |
	x_1^k \rangle $, $A_{n,1} = 0$, and $A_{n,2} = \Lambda(x_1, \dots, x_n) $ are finite-dimensional. This proves the ``only if'' part.
	
	To prove  the ``if'' part, note that each algebra $A_{n,k}$ for $n\ge 2$ and $k\ge 3$ has a quotient algebra isomorphic to $A_{2,3}$. In turn,
	the algebra
	$$
	A_{2,3} =  K\!\langle x,y |\;
	(\alpha x +\beta y)^3 , \alpha, \beta \in K \rangle
	$$
	has a quotient algebra
	$$C = 
	K\!\langle x,\, y\;\;|\;\; x^3,\, y^3,\, x^2 y+xyx+yx^2, \, xy^2+yxy+y^2 x\rangle ,$$
	because 
	$$(\alpha x+\beta y)^3 
	=\alpha^3 x^3 +\beta^3 y^3
	+\alpha^2 \beta (x^2 y+xyx+yx^2)
	+ \alpha\beta^2 (xy^2+yxy+y^2 x)
	$$
	(in fact, the algebras $A_{2,3}$ and $C$ are isomorphic).
	It remains to prove that $C$ is infinite-dimensional. 
	Indeed, 
	one can check that the above
	relations of $C$
	form a Gr\"{o}bner basis
	with respect to the degree-lexicographical order with $x>y>z$. Then the monomials which are not divisible by the leading terms $x^3, x^2y, xy^2, y^3$ of
	the Gr\"{o}bner basis elements are linearly independent in $C$. For instance, the monomials $(xy)^t$ (where $t\ge 1$) are linearly independent,
	so that $C$ is infinite-dimensional.
	\ep\\
	
	\esp As expected, we will give in the sequel a counter-example to the Kurosh problem for Bernstein algebras.
	
	\begin{Example}
		\label{Example: Kurosh for Bernstein}
		 {\rm
			Involving the above  associative algebra $C$ from the proof of Lemma~\ref{Lemma:nil_assoc_algebras}, we  construct the corresponding Bernstein algebra by means of Theorem~\ref{the:Bernstein_via_assoc}. Let $S=K\{x,y\}$ be the two-dimensional vector subspace of $C$ spanned by $x$ and $y$. Let 
			$A = K \times C \times S$ be the Bernstein algebra described in Theorem~\ref{the:Bernstein_via_assoc}.
			Then $A$ is generated by the  set $\{e, (0,x,0), (0,y,0), (0,0,x), (0,0,y)\}$. Now, as in the proof of Theorem \ref{th:Bernstein_via_assoc}, if $v=(0,0,b)\in V$ with $b\in S$, and $u=(0,a,0)\in U$ with $a\in C$, then $L_v^3(u)=(0,ab^3,0))$. But $b$ is a linear combination of $x$ and $y$, so  $b^3=0$. Thus, $L_v^3=0$ for all $v\in V$. Applying Theorem~\ref{th:train_algebra},
			we conclude that $A$ is a train algebra.\\ 
						Let us explicitly write out the train equation for $A$.
			If $t=(\alpha, a, b)\in A$ with $a\in C, b\in S$ and $\alpha = \omega(t)\in K$, then by routine calculations, we have $$t^2-\omega(t)t=(0,2ab, -\omega(t) b), \;\;\;\;\;\; t^3-\omega(t)t^2=t(t^2-\omega(t)t)=(0,2ab^2, 0),$$
			$$(t^3-\omega(t)t^2)(t-\frac 12\omega(t))=t(t^3-\omega(t)t^2)-\frac 12\omega(t)(t^3-\omega(t)t^2)=(0,2\omega(t) ab^3, 0)=(0,0,0).$$
			Hence, $A$ is a train algebra of rank $4$ with train equation $$t^4-\frac 32\omega(t)t^3+\frac 12\omega(t)^2t^2=0.$$
			}\end{Example}
			
	\begin{Note} 
	\label{rem:Kurosh_example_n_k}
	{\rm As far as we know in the literature, the previous algebra $A$ is actually the first example  of an infinite-dimensional finitely generated Bernstein-train algebra.  \\
	More generally, let the associative algebra $A_{n,k}$ of Lemma \ref{Lemma:nil_assoc_algebras} with $n\geq 2$ and $k\geq 3$. Then, by straightforward computations,  Theorem~\ref{the:Bernstein_via_assoc} gives rise to an infinite-dimensional Bernstein-train algebra of rank $k+1$ with $n$ generators.   }
	\end{Note}
	
	\section{Nilpotency in Bernstein algebras and the Jacobian problem}
	
		\label{sec:Jacobian}

	There are various versions of nilpotence in non-associative algebras.
	Let us discuss them in an arbitrary Bernstein algebra $A$. Since $A$ is commutative
	but not necessarily power-associative, the definition of the powers of an element $a\in A$ has several
	variants. But we use here the right principal powers defined inductively by $a^1 = a$ and $a^{n+1} = a^na$.
	
	We say that $a$ is right nilpotent if there
	is an integer $m$ such that $a^m = 0$. The smallest such
	$m$ is called the right nilpotency index of $a$.
	
	We say that $a$ is nilpotent (or, strongly
	nilpotent) if there is an integer $m$ such that each
	product of $m$ copies of $a$, with any arrangement of
	parentheses, is zero. The smallest such $m$ is called the
	nilpotency index of $a$.
	

	Nilpotence conditions in nonassociative algebras are particular cases of similar conditions for multilinear operators. The latter naturally appear in an equivalent formulation of the Jacobian problem due to Yagzhev, see~\cite{Belov_etc} and references therein.

	Let $A$ be a vector space and $\mu : A^d \to A$ be a $d$-linear map, so that the pair $(A,\mu )$ form a (multilinear) algebra (in the sense of general algebra). Let $T^{multilinear}_q (x_1, \dots , x_q)$ be the sum of all multilinear multiple compositions of $\mu $ with $q$ arguments, and let $T_q(x) = T_q^{multilinear}(x,x, \dots, x)$. The operator $\mu $ is called {\em weakly nilpotent} if there exists $q_0>0$ such that  the identity $T_q(x) = 0$ holds for all $x\in A$ and $q\ge q_0$. In this case, $A$ is 
	called\footnote{Note that there is a little inaccuracy in the definition of Yagzhev (or weakly nilpotent) algebra in~\cite[Definition~2.4]{Belov_etc}. Namely, after a verbal definition it is said that ``equation (5) means that $A$ is weakly nilpotent''. In fact, the verbal definition is not equivalent to this condition. In the subsequent theorems in~\cite{Belov_etc}, equation (5) is used as a definition of Yagzhev algebras. The definition given here is equivalent to the equation (5) in~\cite{Belov_etc}.}
	 a Yagzhev algebra~\cite{Belov_etc}.
	
	On the other hand, consider the linear operator $\AdOper_x: y\to \mu (x, \dots, x, y)$ defined by $\mu $ and an element $x\in A$. The algebra $(A,\mu )$ is called {\em Engel } if there exists $p>0$ such that $\AdOper_x$ is nilpotent of degree $p$ for all $x$, 
	that is, $\AdOper_x^p(y) = 0$ for all $x,y\in A$.

	Recall that the famous Jacobian conjecture is equivalent to the following statement:
	
	\begin{Conjecture}[Jacobian conjecture for homogeneous mapping]\!\!\!. \label{conj:Jacobian_hom}
	\\
		Suppose that $F: \Complex^n \to \Complex^n$ is a polynomial automorphism of the form 
		$F = \Identity - H$, where $H$ is a homogeneous automorphism of degree $d\ge 2$.
		If the determinant $j_F$ of the Jacobian matrix $J_F = (\partial F(X)/\partial X)$ of the map $F$ is equal to the constant 1, then the map $F$ has a polynomial inverse. 
	\end{Conjecture}
	
	The conjecture is known to be true for $d=2$~\cite{Wang}, and it is sufficient to prove it for $d=3$~(\cite{Yagzhev}; see also \cite{bass}). 
	
	It follows from the results by Yagzhev that the Jacobian conjecture for homogeneous maps of degree $d$ for polynomials of $n$ variables  is equivalent to the following:
	
	\begin{Conjecture}[Jacobian conjecture in Yagzhev form]\!\!\!.
		\label{conj:Yagzhev}\\
		Suppose that $\mu :A^d\to A$ is a $d$-linear operator on a complex vector space $A$ which is symmetric, that is, for each substitution $\sigma \in S_d$ we have
		$\mu (x_1, \dots , x_d ) = \mu (x_{\sigma 1}, \dots , x_{\sigma d})$.
		Assume also that the system of Capelli identities of order $n+1$ holds in the algebra $(A,\mu )$. If, in addition, the algebra $A$ is Engel, then it is Yagzhev, that is, 
		$\mu $ is weakly nilpotent. 
	\end{Conjecture}
	
	In the notation of Conjecture~\ref{conj:Jacobian_hom}, here $\mu $ is a linearization of $H$ so that $ H(X) = \mu (X, \dots, X)$.

If we omit the condition about the Capelli identities, we get a natural version of the Jacobian conjecture  for infinite number $n$ of variables. For quadratic algebras, we get the following conjecture~\cite[2.2.2]{Belov_etc}. This conjecture is a natural generalization of the Jacobian conjecture for quadratic mappings proved by Wang~\cite{Wang}.

\begin{Conjecture}[Generalized Jacobian conjecture for quadratic mappings]\!\!\!.\\
	Let $A$ be a commutative binary complex algebra. If $A$ is Engel, then  $A$ is Yagzhev.
\end{Conjecture}

For example, suppose that $A$ is power-associative (e.g., a Jordan algebra).
Then we have $T_n(x) = C_n x^n = C_n {\AdOper}_x^{n-1} (x)$, where $C_n>0$. If $A$ is Engel, then 
${\AdOper}_x^{n-1} (x) =0$ for large enough $n$, so that $T_n(x)=0$ and $A$ is Yagzhev. So, the generalized 
Jacobian conjecture for quadratic mappings holds for $A$.\\ 

Consider now the generalized Jacobian conjecture for quadratic mappings in two classes of algebras, namely, commutative algebras with the identity $(x^2)^2=0$  and subalgebras of Bernstein algebras. Recall that the previous identity holds 
in the barideal of a Bernstein algebra.

\begin{Theorem}
\label{th:Bernstein_Engel_Yagzhev}
Let $K$ be an arbitrary field of characteristic different from $2$. 
Assume that the identity $(x^2)^2=0$ holds in a commutative algebra $A$ over $K$.
Then the following conditions are equivalent:

(i) $A$ is a nil-algebra of bounded nil-index;

(ii) $A$ is Engel;

(iii) $A$ Yagzhev.
\end{Theorem}

In particular, one concludes that the generalized Jacobian conjecture for quadratic mappings holds in $A$
with no restriction for the field characteristic, except for the above.

To prove the theorem, we need 
	
	\begin{Lemma}
		\label{lem:nilpotency}
		Let $N$ be a commutative algebra over a field $K$ of characteristic different from $2$. 
		Assume that $N$ satisfies the identity $(x^2)^2=0$ (this is in particular the case when $N$ is the barideal of a Bernstein algebra). Let $a\in N$.
	 Then
		\begin{enumerate}
			\item With the exception of the 
			principal powers
			$a^m$ of $a$, every product of $m\geq 4$ copies of $a$, with
			any arrangement of parentheses, is zero.
			\item If $a$ is right nilpotent of right nilpotency index
			$m \geq  2$, then $a$ is nilpotent of the same nilpotency
			index $m$.
		\end{enumerate}
	\end{Lemma}
	
	\pr 
     Let us first prove that 
  \begin{equation}
\label{eq:x^ix^j=0}
x^ix^j=0 \mbox{ for all } x\in N \mbox{ and } i,j\geq 2.
\end{equation}

Indeed, a partial linearization of the identity $(x^2)^2=0$ yields the identities  
\begin{equation}
\label{eq:x^2(xy)=0}
x^2(xy)=0,
\end{equation}
\begin{equation}
\label{eq:2(xy)(xz)+x^2(yz)=0}
2(xy)(xz)+x^2(yz)=0.
\end{equation}
Putting $y=x^{j-1}$ in (\ref{eq:x^2(xy)=0}), we get 
$$x^2x^j=0$$
for each $j\ge 2$. Then taking $y=x^{i-1}, z= x^{j-1}$ in~(\ref{eq:2(xy)(xz)+x^2(yz)=0}), 
we obtain by induction $x^ix^j =0$ for each $i,j\ge 3$.\\
	A product of $m$ copies of $a$, with a chosen distribution of parentheses, can be considered as a
	nonassociative monomial $M_m(a)$ in $a$ of degree $m$.   \\
		1. For $m\leq 3,  M_m(a)$ is a principal power. The assertion is obvious for $m=4$, because $(a^2)^2=0$. 
		Let $m \geq 5$ and assume by induction that every monomial $M_t(a)\ne a^t$ in $a$ of degree $t< m$ vanishes. 
		Let $M_{m}(a)$ be a monomial in $a$ of degree $m$.
		Write the representation of the monomial $M_{m}(a)$ as
		the product $M_{m}(a) = M_i(a)  M_j(a)$ of two monomials 		$M_i(a)$ and $M_j(a)$ in $a$ of lesser degrees $i$ and $j$.
		By the induction hypothesis, if $M_i(a)$ and $M_j(a)$ are nonzero, then they are 
		right principal powers. In this case, when $i=1$ or $j=1$, then 
		$M_{m}(a) = a^{m-1} a=a^m$ is also a right principal power. Otherwise,
		$M_{m}(a) = a^i a^j$ for some $i,j\ge 2$, so that $M_{m}(a)=0$ by (\ref{eq:x^ix^j=0}).
		
		2. It follows immediately from the first assertion, since  the only monomial in $a$ of degree $m$, which may be different from zero,
		is precisely the right principal power $a^m$.\ep\\


{\it Proof of Theorem~{\rm \ref{th:Bernstein_Engel_Yagzhev}}}: 
By Bayara's
Lemma~\ref{lem:Bayara_on_nil}, the conditions $(i)$ and $(ii)$ are equivalent. Moreover, it follows from Lemma~\ref{lem:nilpotency} that for each $x\in A$ and $q\ge 2$, we have
$T_q(x) = 2^{q-2} x^q$. Since the characteristic of the ground field $K$ is not 2, 
the identities $T_q(x)=0$ and $x^q =0$ are equivalent. This proves the equivalence $(i)\Longleftrightarrow (iii)$.
\ep

\begin{Corollary}
\label{cor:jacobian-bernstein}
	Suppose that $A$ is a subalgebra of a Bernstein algebra (over an arbitrary field of characteristic different from $2$).  
	Then the  generalized Jacobian conjecture for quadratic mappings holds for $A$.  
\end{Corollary}

\pr
Assume that $A$ is a subalgebra of a Bernstein algebra $B$. Let $x\in A$. If $A$ is 
Engel, then $x^n = 0$ for some $n>0$. So, $w(x^n) = w(x)^n = 0$, and $x\in \ker w = N(B)$. Hence $A$ is a subalgebra of $N(B)$. Since the identity  $(x^2)^2=0$ holds in the barideal
$N(B)$, we apply the implication $(ii)\Longrightarrow (iii)$ of Theorem~\ref{th:Bernstein_Engel_Yagzhev}
to conclude that $A$ is Yagzhev. 
\ep\\

\vspace*{0,5cm}

{\large \bf Acknowledgment:}\\
The second author wishes to thank Professor Moussa Ouattara for a number of useful discussions and for reading an earlier
draft of this article. 

The work of D. Piontkovski has been supported by the grant of the Russian Science Foundation, RSF 22-21-00912.

\end{document}